# Tetraplectic structures compatible with local quaternionic toric actions

Panagiotis Batakidis[†]  Ioannis Gkeneralis[†]
batakidis@math.auth.gr  igkeneralis@math.auth.gr

[†] Department of Mathematics, Aristotle University of Thessaloniki

**Abstract**

This paper introduces a quaternionic analogue of toric geometry by developing the theory of local $Q^n := \text{Sp}(1)^n$-actions on $4n$-dimensional manifolds, modeled on the regular representation. We identify obstructions that measure the failure of local properties to globalize and define two invariants: a combinatorial invariant called the characteristic pair and a cohomological invariant called the Euler class, which together classify local quaternionic torus actions up to homeomorphism. We also study tetraplectic structures in quaternionic toric geometry by introducing locally generalized Lagrangian-type toric fibrations and show that such fibrations are locally modeled on $\mathbb{R}^n \times Q^n$ using a quaternionic version of the Arnold–Liouville theorem. In the last part, we show that orbit spaces of these actions acquire the structure of quaternionic integral affine manifolds with corners and Lagrangian overlaps, and we classify such spaces by establishing a quaternionic Delzant-type theorem.

# Contents

# 1   Introduction

Toric geometry lies at the intersection of algebraic geometry, symplectic geometry, and combinatorics. In its classical formulation, it concerns the study of algebraic varieties or symplectic manifolds equipped with effective actions of a compact torus $T^n$. The geometry and topology of such spaces are deeply connected to convex polytopes via moment maps, which reflect the orbit structure and symmetries of the action (see [Buchstaber and Panov, 2014], [Buchstaber and Panov, 2002], [Buchstaber and Panov, 2000]).

Two major perspectives shape the field. The first, arising from algebraic geometry, interprets toric varieties as varieties associated with fans—collections of strongly convex rational polyhedral cones. The second, rooted in symplectic geometry, views toric manifolds as symplectic manifolds with effective Hamiltonian $T^n$-actions, where the image of the moment map is a convex polytope. Delzant's theorem [Delzant, 1988], together with the convexity results of Atiyah [Atiyah, 1982] and Guillemin–Sternberg [Guillemin and Sternberg, 1982], establishes a complete classification: up to equivariant symplectomorphisms, there is a bijection between compact symplectic toric manifolds and Delzant polytopes.

In [Davis and Januszkiewicz, 1991], Davis and Januszkiewicz introduced a topological analogue of toric varieties, namely *quasitoric manifolds*. Quasitoric manifolds are smooth manifolds equipped with locally standard $T^n$-actions (left multiplication) whose orbit spaces are simple convex polytopes. After their work, more topological generalizations appeared in the literature (see [Buchstaber and Panov, 2000, Buchstaber and Panov, 2002] for further background).

In this paper, we generalize the ideas of [Yoshida, 2011] in the quaternionic setting. We first introduce the notion of a *local quaternionic torus action* modeled on the regular representation of the compact Lie group $Q := \mathrm{Sp}(1) \cong S^3$ (we call $Q$ the *quaternionic tori*). This generalizes the concept of a *locally regular quaternionic torus action*, as studied in the works of Scott on quaternionic toric varieties [Scott, 1995] and Hopkinson's development of *quoric manifolds* [Hopkinson, 2012].

More precisely, our starting point is the introduction of weakly regular atlases on $4n$-dimensional manifolds, whose transition functions encode the local symmetry data via the automorphism group $\mathrm{Aut}(Q^n)$. This setup allows us to define an orbit space $B_M$, which inherits a natural structure of an $n$-dimensional topological manifold with corners. The obstruction to globalizing such a local action is captured by a Čech cohomology class in $H^1(B_M; \mathrm{Aut}(Q^n))$, generalizing similar phenomena in toric topology. These orbit spaces provide the natural base spaces for our subsequent study of generalized quaternionic toric fibrations.

To get a topological classification of local quaternionic torus actions, we introduce two fundamental invariants associated with such actions. The first invariant is a *characteristic pair* $(P, \mathcal{L})$, defined via bundle-theoretic data on the orbit space $B_M$. It encodes the isotropy structure of the action in terms of a principal $\mathrm{Aut}(Q^n)$-bundle and a system of sublattices over the boundary strata. The second invariant is an *Euler class*, that is a Čech cohomology class with values in a sheaf of sections of a canonical torus bundle. These invariants generalize the role of characteristic functions and Euler classes in the classification of torus manifolds, that appeared in [Yoshida, 2011]. Our first main result shows that local quaternionic torus actions are completely determined by their characteristic pair and Euler class:

**Theorem** (Main Theorem 1). *(Classification of local quaternionic torus actions) Two compact $4n$-dimensional manifolds equipped with local quaternionic torus actions are homeomorphic if and only if their characteristic pairs are isomorphic and their Euler classes coincide under the induced pullback.*

An important subclass of such actions arises when the underlying manifold is equipped with a *tetraplectic structure*, that is a closed, nondegenerate 4-form that is invariant under the action. This structure, introduced by Foth [Foth, 2002], serves as the quaternionic analogue of a symplectic form, and allows for the definition of *tri-moment maps*, which generalize the classical moment maps of Hamiltonian torus actions. Our focus is on understanding how these local actions, when combined with tetraplectic geometry, give rise to generalized toric-like fibrations and a new class of topological and geometric models.

A key analytical ingredient in this theory is a quaternionic generalization of the Arnold–Liouville theorem. After we define and study *locally generalized Lagrangian-type toric fibrations*, we prove that such fibrations are locally tetraplectomorphic to a standard product model

$$(\mathbb{R}^n \times Q^n, \psi_{\mathbb{R}^n \times Q^n}),$$

where $\psi$ is a canonical closed 4-form encoding the tetraplectic structure. This result mirrors the role of action-angle coordinates in classical integrable systems.

We then establish a necessary and sufficient condition for a tetraplectic $Q^n$-manifold to admit such a fibration: the base must carry a quaternionic integral affine structure, and the transition data must be encoded by a class in the Čech cohomology of a sheaf of Lagrangian sections. This cohomological viewpoint leads naturally to a quaternionic analogue of the Delzant classification theorem. Specifically, we classify such spaces in terms of their orbit space and Lagrangian class, recovering and extending examples such as the quaternionic projective spaces studied in [Gentili et al., 2019]. Our second main result is the following.

**Theorem** (Main Theorem 2). *A $4n$-dimensional manifold equipped with a local quaternionic torus action admits a tetraplectic structure if and only if the orbit space carries a quaternionic integral affine structure and the overlap data defines a Lagrangian Čech cocycle.*

It should be noted that analogous results have been established in the classical (complex) toric setting, as seen in the works of Duistermaat [Duistermaat, 1980], Mishachev [Mishachev, 1996], Symington [Symington, 2003], Gay–Symington [Gay and Symington, 2009], and Yoshida [Yoshida, 2011].

The paper is organized as follows. In Section 2, we recall the definition of the regular representation and introduce the notion of a regular atlas. We then generalize this to weakly regular atlases, which define local $Q^n$-actions modeled on the regular representation. We also examine the orbit space and orbit maps associated with local $Q^n$-actions and present an obstruction to their globalization via Čech cohomology with values in $Aut(Q^n)$. Section 3 is devoted to the construction of the characteristic pair and the associated canonical model. In Section 4, we give a topological criterion for the existence of global sections of the orbit map and establish a classification theorem for manifolds equipped with local $Q^n$-actions. In the final Section, we turn to the subclass of such manifolds that admit tetraplectic structures. Here, we prove two important results: a realization theorem and a classification theorem in the quaternionic-tetraplectic setting.

**Acknowledgments:** We would like to thank Nikolas Adaloglou and Efstratios Prassidis for useful discussions on the subject.

## 2 Local quaternionic torus actions

### 2.1 The regular representation

Let $Q := S^3 \subset \mathbb{H}$ be the group of unit quaternions and $Q^n$ be the n-dimensional quaternionic torus which is considered acting on $\mathbb{H}^n$ by coordinatewise quaterninonic multiplication,
$$((s_1, \ldots, s_n), (h_1, \ldots, h_n)) \mapsto (s_1 h_1, \ldots, s_n h_n).$$
This action will be called the *regular $Q^n$-action on the regular n-corner* $\mathbb{H}^n$, also known as the *regular representation* of $Q^n$. Regular representations are characterized by the following general properties:

1. The orbits are products of 3-spheres and the origin is the only fixed point.

2. The orbit space $\mathbb{H}^n/Q^n$ is the positive cone $\mathbb{R}^n_{\geq 0} = \{(h_1, \ldots, h_n) \mid h_i \geq 0\}$, with radii $|h|$, $h \in \mathbb{H}^n$ and the quotient map is given by
$$(h_1, \ldots, h_n) \mapsto (|h_1|^2, \ldots, |h_n|^2).$$

3. The conjugacy classes are naturally associated to the faces of the orbits. We refer to them as *isotropy classes*.

### 2.2 Locally regular torus actions

**Definition 2.1.** Let $M^{4n}$ be a 4n-dimensional smooth manifold equipped with a smooth $Q^n$-action. A *regular coordinate neighborhood* of $M^{4n}$ is a triple $(U, \rho, \varphi)$ consisting of a $Q^n$-invariant open set $U$ of $M^{4n}$, an automorphism $\rho$ of $Q^n$, and a $\rho$-equivariant diffeomorphism $\varphi$ from $U$ to some $Q^n$-invariant open subset of $\mathbb{H}^n$. The action of $Q^n$ on $M^{4n}$ is said to be *locally regular* if every point in $M^{4n}$ lies in some regular coordinate neighborhood, and an atlas which consists of regular coordinate neighborhoods is called *regular atlas*.

**Example 2.2.** A quoric manifold is a smooth manifold equipped with a locally regular torus action whose orbit space is combinatorially isomorphic (as a manifold with corners) to a simple polytope. That is, there exists a face-preserving homeomorphism between the orbit space and the polytope. Quoric manifolds were introduced by Hopkinson in [Hopkinson, 2012], and generalize topologically quasitoric manifolds (for the later see [Davis and Januszkiewicz, 1991]).

For a topological space which is not necessarily equipped with a global $Q^n$-action we consider a generalization of a regular atlas.

**Definition 2.3.** Let $M$ be a compact Hausdorff space. A *weakly regular ($0 \leq r \leq \infty$) atlas* of $M$ is an atlas $\{(U_\alpha^M, \varphi_\alpha^M)\}_{\alpha \in \mathcal{A}}$ which satisfies the following properties:

1. for each $\alpha$, $\varphi_\alpha^M$ is a homeomorphism from $U_\alpha^M$ to an open set of $\mathbb{H}^n$ which is invariant with respect to the regular representation of $Q^n$.

2. for each nonempty overlap $U_{\alpha\beta}^M := U_\alpha^M \cap U_\beta^M$

    a. $\varphi_\alpha^M(U_\alpha^M)$ and $\varphi_\beta^M(U_\beta^M)$ are also invariant with respect to the regular representation of $Q^n$.

    b. there exists an automorphism $\rho_{\alpha\beta}$ of $Q^n$ such that the overlap map $\varphi_{\alpha\beta}^M := \varphi_\alpha^M \circ (\varphi_\beta^M)^{-1}$ is $C^r$, is a $\rho_{\alpha\beta}$-equivariant diffeomorphism of $\mathbb{H}^n$, where the $Q^n$ acts by the regular representation.

Two weakly regular atlases $\{(U_\alpha^M, \varphi_\alpha^M)\}_{\alpha \in A}$ and $\{(V_\beta^M, \psi_\beta^M)\}_{\beta \in B}$ of $M$ are said to be *equivalent* if on each nonempty overlap $U_\alpha^M \cap V_\beta^M$, there exists an automorphism $\rho$ of $Q^n$ such that $\varphi_\alpha^M \circ (\psi_\beta^M)^{-1}$ is a $\rho$-equivariant diffeomorphism. We call an equivalence class of weakly regular atlases a *local $Q^n$-action on $M$ modeled on the regular representation*. If no confusion is possible, we will call it *local torus/$Q^n$-action on $M$* and denote it by $\mathcal{Q}$.

For each local torus action $\mathcal{Q}$ there uniquely exists a weakly regular atlas in $\mathcal{Q}$ that contains all weakly regular atlases in $\mathcal{Q}$. We call it a *maximal* weakly regular atlas in $\mathcal{Q}$.

**Definition 2.4.** Let $(M_i, \mathcal{Q}_i)$ $(i = 1, 2)$ be a 4n-dimensional manifold equipped with a local $Q^n$-action $\mathcal{Q}_i$ and $\{(U_\alpha^{M_1}, \varphi_\alpha^{M_1})\}_{\alpha \in A} \in \mathcal{Q}_1$ and $\{(U_\beta^{M_2}, \varphi_\beta^{M_2})\}_{\beta \in B} \in \mathcal{Q}_2$ be the maximal weakly regular atlases of $M_1$ and $M_2$. A $C^r$ *isomorphism* $f_M : (M_1, \mathcal{Q}_1) \to (M_2, \mathcal{Q}_2)$ between $(M_1, \mathcal{Q}_1)$ and $(M_2, \mathcal{Q}_2)$ is a diffeomorphism $f_M : M_1 \to M_2$ which satisfies the following conditions: for each nonempty overlap $U_\alpha^{M_1} \cap f_M^{-1}(U_\beta^{M_2})$,

(i) $\varphi_\alpha^{M_1}(U_\alpha^{M_1} \cap f_M^{-1}(U_\beta^{M_2}))$ and $\varphi_\beta^{M_2}(f_M(U_\alpha^{M_1}) \cap U_\beta^{M_2})$ are invariant with respect to the regular representation of $Q^n$ and

(ii) there exists an automorphism $\rho$ of $Q^n$ such that $\varphi_\beta^{M_2} \circ f_M \circ (\varphi_\alpha^{M_1})^{-1}$ is $\rho$-equivariant.

$(M_1, \mathcal{Q}_1)$ and $(M_2, \mathcal{Q}_2)$ are said to be $C^r$ *isomorphic* if there exists a $C^r$ isomorphism between $(M_1, \mathcal{Q}_1)$ and $(M_2, \mathcal{Q}_2)$.

**Remark 2.5.** We will see in the next subsection (Proposition 2.12) that the existence of a local $Q^n$-action does not necessarily imply an existence of a global $Q^n$-action. The obstruction described by a Čech cohomology class determined by automorphisms $\rho_{\alpha\beta}$ of $Aut(Q^n)$.

## 2.3 Orbit spaces and orbit maps

Let $(M^{4n}, \mathcal{Q})$ be a 4n-dimensional manifold equipped with a local $Q^n$-action. For $(M^{4n}, \mathcal{Q})$ we define the orbit space and the orbit map in the following way. Let $\{(U_\alpha^M, \varphi_\alpha^M)\}_{\alpha \in \mathcal{A}} \in \mathcal{Q}$ be the maximal weakly regular atlas of $M^{4n}$.

We endow each quotient space $\varphi_\alpha^M(U_\alpha^M)/Q^n$ with the quotient topology induced from the topology of $\varphi_\alpha^M(U_\alpha^M)$ by the natural projection

$$\pi_\alpha : \varphi_\alpha^M(U_\alpha^M) \to \varphi_\alpha^M(U_\alpha^M)/Q^n.$$

By Definition 2.3, for each overlap $U_{\alpha\beta}^M$, $\varphi_{\alpha\beta}^M$ induces a homeomorphism from $\varphi_\beta^M(U_{\alpha\beta}^M)/Q^n$ to $\varphi_\alpha^M(U_{\alpha\beta}^M)/Q^n$.

Define two elements $b_\alpha \in \varphi_\alpha^M(U_\alpha^M)/Q^n$ and $b_\beta \in \varphi_\beta^M(U_\beta^M)/Q^n$ to be equivalent, and write $b_\alpha \sim_{\text{orb}} b_\beta$, if

$$b_\alpha \in \varphi_\alpha^M(U_{\alpha\beta}^M)/Q^n, \quad b_\beta \in \varphi_\beta^M(U_{\alpha\beta}^M)/Q^n$$

and the homeomorphism induced by $\varphi_{\alpha\beta}^M$ sends $b_\beta$ to $b_\alpha$. It is an equivalence relation in the disjoint union

$$\coprod_\alpha (\varphi_\alpha^M(U_\alpha^M)/Q^n).$$

We define the *orbit space* $B_M$ of the local $Q^n$ action to be the quotient space

$$B_M := \coprod_\alpha (\varphi_\alpha^M(U_\alpha^M)/Q^n)/\sim_{orb}$$

together with the quotient topology. It is easy to see that $B_M$ is a Hausdorff space and $\{\varphi_\alpha^M(U_\alpha^M)/Q^n\}_{\alpha \in \mathcal{A}}$ is an open covering of $B_M$. By construction of $B_M$, the map

$$\coprod_\alpha \pi_\alpha \circ \varphi_\alpha^M : \coprod_\alpha U_\alpha^M \to \coprod_\alpha (\varphi_\alpha^M(U_\alpha^M)/Q^n)$$

induces an open map from $M^{4n}$ to $B_M$. We call it the *orbit map* of the local $Q^n$-action and denote it by
$$\pi_M : M^{4n} \to B_M$$

**Definition 2.6** (see [Davis, 2008]). A space $X$ is an *n-manifold with corners* if it is Hausdorff, second countable space equipped with an atlas of open sets homeomorphic to open subsets of $\mathbb{R}^n_+$ such that the overlap maps are local homeomorphisms that preserve the natural stratification of $\mathbb{R}^n_+$.

**Proposition 2.7.** $B_M$ *is endowed with a structure of an n-dimensional manifold with corners.*

*Proof.* A structure of an $n$-dimensional manifold with corners on $B_M$ is constructed as follows. We put
$$U_\alpha^B := \varphi_\alpha^M(U_\alpha^M)/Q^n.$$
The restriction of $\pi_{\mathbb{H}^n}$ to $\varphi_\alpha^M(U_\alpha^M)$ induces a homeomorphism from $U_\alpha^B$ to the open subset $\pi_{\mathbb{H}^n}(\varphi_\alpha^M(U_\alpha^M))$ of $\mathbb{R}^n_{\geq 0}$, which is denoted by $\varphi_\alpha^B$. By construction, on each overlap $U_{\alpha\beta}^B := U_\alpha^B \cap U_\beta^B$, the overlap map
$$\varphi_{\alpha\beta}^B := \varphi_\alpha^B \circ (\varphi_\beta^B)^{-1} : \pi_{\mathbb{H}^n}(\varphi_\beta^M(U_{\alpha\beta}^M)) \to \pi_{\mathbb{H}^n}(\varphi_\alpha^M(U_{\alpha\beta}^M))$$
preserves the natural stratifications of $\pi_{\mathbb{H}^n}(\varphi_\alpha^M(U_{\alpha\beta}^M))$ and $\pi_{\mathbb{H}^n}(\varphi_\beta^M(U_{\alpha\beta}^M))$. Thus, $\{(U_\alpha^B, \varphi_\alpha^B)\}_{\alpha \in \mathcal{A}}$ is the desired atlas. □

The fact that the orbit space $B_M$ is a smooth n-manifold with corners implies that it carries a natural orbit stratification. We describe this stratification on $B_M$ as follows. Let $\mathcal{S}^{(k)}B_M$ denote the codimension-$k$ stratum of the stratification. The set of points where exactly $k$ local coordinates vanish forms a smooth submanifold $\mathcal{S}^{(k)}B \subset B_M$ of codimension $k$. The decomposition
$$B_M = \bigsqcup_{k=0}^{n} \mathcal{S}^{(k)}B_M$$
defines a natural stratification of $B_M$, where each stratum $\mathcal{S}^{(k)}B_M$ corresponds to the set of points lying on codimension-$k$ faces of $B_M$. This stratification induces a canonical filtration
$$\mathcal{S}^{(0)}B_M \subset \mathcal{S}^{(1)}B_M \subset \cdots \subset \mathcal{S}^{(n)}B_M = B_M,$$
where $\mathcal{S}^{(k)}B_M := \bigcup_{j=0}^{k} S_j$ denotes the union of all strata of codimension less than or equal to $k$. Our treatment of this stratification is based on the classical bundle theory; see [Lawson and Michelsohn, 1989]. This particular stratification will be needed in a later section to define characteristc pairs and the canonical model.

**Remark 2.8.** In the case of a locally regular quaternionic torus action, $\mathcal{S}^{(0)}B_M$ corresponds to the fixed point set.

**Remark 2.9.** It would be interesting to develop an approach to local quaternionic torus actions based on Quinn's theory of homotopically stratified spaces and controlled topology, rather than relying solely on bundle-theoretic techniques. The defining conditions of Quinn's stratified spaces are homotopy-theoretical, from which geometric-topological properties can be deduced. Such a framework could possibly provide finer tools for analyzing the local structure near singular strata in the orbit space.

**Remark 2.10.** The atlas $\{(U_\alpha^B, \varphi_\alpha^B)\}_{\alpha \in \mathcal{A}}$ of $B_M$ constructed in Proposition 2.7, has the following property: for each $\alpha$, $U_\alpha^M = \pi_M^{-1}(U_\alpha^B)$, $\varphi_\alpha^M(U_\alpha^M) = \pi_{\mathbb{H}^n}^{-1}(\varphi_\alpha^B(U_\alpha^B))$ and the following diagram commutes

$$\begin{array}{ccccccc}
M^{4n} & \supset & \pi_M^{-1}(U_\alpha^B) & \xrightarrow{\varphi_\alpha^M} & \pi_{\mathbb{H}^n}^{-1}(\varphi_\alpha^B(U_\alpha^B)) & \subset & \mathbb{H}^n \\
\downarrow \pi_M & & \downarrow \pi_M & & \downarrow \pi_{\mathbb{H}^n} & & \downarrow \pi_{\mathbb{H}^n} \\
B_X & \supset & U_\alpha^B & \xrightarrow{\varphi_\alpha^B} & \varphi_\alpha^B(U_\alpha^B) & \subset & \mathbb{R}^n_{\geq 0}
\end{array}$$

**Remark 2.11.** Let $(M_1^{4n}, \mathcal{Q}_1)$ and $(M_2^{4n}, \mathcal{Q}_2)$ be $4n$-dimensional manifolds equipped with local $Q^n$-actions. Consider the isomorphism defined in 2.4,

$$f_M : (M_1^{4n}, \mathcal{Q}_1) \to (M_2^{4n}, \mathcal{Q}_2).$$

Then, $f_M$ induces a stratification preserving homeomorphism $f_B : B_{M_1} \to B_{M_2}$ such that

$$f_M = f_B \circ \pi_{M_1}.$$

Later, we will classify local quaternionic torus actions up to homeomorphisms.

Let $\{(U_\alpha^B, \varphi_\alpha^B)\}_{\alpha \in \mathcal{A}}$ the atlas of $B_M$. The automorphisms $\rho_{\alpha\beta}$ of Definition 2.3 form a Čech 1-cocycle $\{\rho_{\alpha\beta}\}$ on $\mathcal{U} := \{U_\alpha^B\}_{\alpha \in \mathcal{A}}$ with values in $Aut(Q^n)$. Hence $\{\rho_{\alpha\beta}\}$ determines a cohomology class in the first Čech cohomology $H^1(B_M, Aut(Q^n))$.

**Proposition 2.12.** *A local $Q^n$-action $\mathcal{Q}$ on $M^{4n}$ is induced by some locally regular $Q^n$-action if and only if $\{\rho_{\alpha\beta}\}$ is cohomologous to the trivial Čech 1-cocycle in $H^1(B_M; Aut(Q^n))$, where the trivial Čech 1-cocycle is the one whose values on all open sets are equal to the identity map of $Q^n$.*

*Proof.* Suppose first that the local $Q^n$-action $\mathcal{Q}$ on $M^{4n}$ is induced by a globally defined locally regular $Q^n$-action. Then there exists a locally regular atlas $\{(U_\alpha^M, \rho_\alpha, \varphi_\alpha^M)\}$ in $\mathcal{Q}$ such that the overlap maps $\varphi_{\alpha\beta}^M := \varphi_\alpha^M \circ (\varphi_\beta^M)^{-1}$ are $\rho_{\alpha\beta}$-equivariant for $\rho_{\alpha\beta} \in \mathrm{Aut}(Q^n)$ of the form

$$\rho_{\alpha\beta} = \rho_\alpha \circ \rho_\beta^{-1}$$

for some $\rho_\alpha, \rho_\beta \in \mathrm{Aut}(Q^n)$. This implies that $\{\rho_{\alpha\beta}\}$ is cohomologous to the trivial Čech 1-cocycle in $H^1(B_M; \mathrm{Aut}(Q^n))$ (for details see [Lawson and Michelsohn, 1989]).

Conversely, assume that the cocycle $\{\rho_{\alpha\beta}\}$ is cohomologous to the trivial cocycle in $H^1(B_M; \mathrm{Aut}(Q^n))$. Then there exists a Čech 0-cochain $\{\rho_\alpha\}$ such that

$$\rho_{\alpha\beta} = \rho_\alpha \circ \rho_\beta^{-1}.$$

Define new charts $\tilde{\varphi}_\alpha^M := \varphi_\alpha^M \circ \rho_\alpha^{-1}$ on $U_\alpha^M$. On overlaps $U_{\alpha\beta}^M$, we compute:

$$\tilde{\varphi}_\alpha^M \circ (\tilde{\varphi}_\beta^M)^{-1} = \varphi_\alpha^M \circ \rho_\alpha^{-1} \circ \rho_\beta \circ (\varphi_\beta^M)^{-1} = \varphi_\alpha^M \circ (\varphi_\beta^M)^{-1} \circ (\rho_\beta^{-1} \circ \rho_\alpha) = \rho_{\alpha\beta} \circ (\rho_\beta^{-1} \circ \rho_\alpha) = \mathrm{id}.$$

Therefore, the new transition maps are trivial and the action becomes locally regular in the sense that the transition automorphisms are globally trivial. Therefore, the local $Q^n$-action is induced by a globally defined locally regular $Q^n$-action. □

## 3 The canonical model

In this section, we introduce characteristic pairs and construct a canonical model from characteristic pairs. We also show that the characteristic pair is associated to a local $Q^n$-action.

### 3.1 The integral lattice

The Lie algebra $\mathfrak{q} = \mathfrak{su}(2)$, consists of $2 \times 2$ traceless anti-Hermitian matrices and is spanned by the generators

$$T_1 = \frac{i}{2}\sigma_1, \quad T_2 = \frac{i}{2}\sigma_2, \quad T_3 = \frac{i}{2}\sigma_3$$

where $\sigma_i$ are the Pauli matrices. The Lie bracket is

$$[T_i, T_j] = \epsilon_{ijk} T_k,$$

and for any element $X = \alpha_1 T_1 + \alpha_2 T_2 + \alpha_3 T_3$, $\alpha_i \in \mathbb{R}$, the exponential map is denoted by $\exp(X) = e^X$. Since $X$ is traceless and anti-Hermitian, we can write it in the form

$$X = i\theta \hat{n} \frac{\vec{\sigma}}{2},$$

where $\theta = \sqrt{\alpha_1^2 + \alpha_2^2 + \alpha_3^2}$ and $\hat{n} = (n_1, n_2, n_3)/\theta$ is a unit vector.

Using the well-known identity for the exponential of a traceless anti-Hermitian matrix expressed via Pauli matrices, we obtain

$$e^X = e^{i\theta \hat{n} \cdot \vec{\sigma}/2} = \cos\left(\frac{\theta}{2}\right) I + i \sin\left(\frac{\theta}{2}\right) (\hat{n} \cdot \vec{\sigma}).$$

This describes a general element of $SU(2)$ as a rotation in quaternionic form.

The matrix $e^X$ is the identity if and only if

$$cos(\frac{\theta}{2}) = 1 \quad and \quad sin(\frac{\theta}{2}) = 0.$$

Thus, an $X \in \mathfrak{q}$ is in the kernel of exp if

$$\theta = 4k\pi, \quad k \in \mathbb{Z}.$$

Therefore, the integral lattice is

$$\Lambda = 4\pi\mathbb{Z} \cdot T_3.$$

We know [Hopkinson, 2012, Proposition 2.1.21] that $\text{Aut}(Q^n)$ consists of inner automorphisms (conjugation), Weyl group actions (sign change) and symmetric group actions (permutation). Recall that $SU(2)$ has no outer automorphisms. Hence, up to inner automorphisms, we have that

$$\text{Aut}(Q^n) \cong S_n \ltimes (\mathbb{Z}/2\mathbb{Z})^n$$

where $(\mathbb{Z}/2\mathbb{Z})^n$ is the Weyl group of $SU(2)$, acting by sign "flips", e.g. $T_3 \mapsto -T_3$ and $S_n$ is the symmetric group. Therefore, we can define the integral lattice of the Lie algebra $\mathfrak{q}^n$ to be

$$\Lambda_Q = (4\pi\mathbb{Z})^n = \{(4\pi k_1, ..., 4\pi k_n), k_i \in \mathbb{Z}.\}$$

The Weyl group $(\mathbb{Z}/2\mathbb{Z})^n$ acts on each coordinate, preserving the integral lattice. Permutation of the coordinates via $S_n$ also preserves the lattice. Together, they form the natural lattice automorphism group,

$$GL(\Lambda_Q) \cong S_n \ltimes (\mathbb{Z}/2\mathbb{Z})^n.$$

Consequently, we obtain

$$GL(\Lambda_Q) \cong \text{Aut}(Q^n).$$

**Remark 3.1.** Inner automorphisms act trivially on the lattice: conjugation preserves the center, which is where lattice elements lie via the exponential map.

## 3.2 Characteristic pairs

In this section, we construct a *characteristic pair* associated to a local $Q^n$-action on a manifold. This consists of a principal $\text{Aut}(Q^n)$-bundle and a rank-one sublattice bundle over the codimension-one stratum of the orbit space. The general bundle-theoretic formalism used here is adapted from [Yoshida, 2011, Section 4.1].

Let $(M^{4n}, \mathcal{Q})$ be a 4n-dimensional manifold equipped with a local $Q^n$-action. Let $\{(U_\alpha^M, \varphi_\alpha^M)\}_{\alpha \in \mathcal{A}} \in \mathcal{Q}$ denote the maximal weakly regular atlas. It induces an atlas $\{(U_\alpha^B, \varphi_\alpha^B)\}_{\alpha \in \mathcal{A}}$ of the orbit space $B_M$, satisfying the compatibility property described in Remark 2.10.

The transition data $\{\rho_{\alpha\beta}\}$ on overlaps of this atlas form a Čech 1-cocycle with values in $Aut(Q^n)$, thereby defining a principal $Aut(Q^n)$-bundle over $B_M$:

$$\pi_{P_M} : P_M \to B_M,$$

explicitly constructed as

$$P_M := \left(\bigsqcup_\alpha U_\alpha^B \times Aut(Q^n)\right) / \sim_P, \tag{3.1}$$

where the equivalence relation $\sim_P$ is given by

$$(b_\alpha, h_\alpha) \sim_P (b_\beta, h_\beta) \quad \text{if } b_\alpha = b_\beta \in U_{\alpha\beta}^B, \ h_\alpha = \rho_{\alpha\beta} \circ h_\beta.$$

The bundle projection $\pi_{P_M}$ is defined in the obvious way. For each $\alpha$, every point in $\pi_{P_M}^{-1}(U_\alpha^B)$ has a unique representative which lies in $U_\alpha^B \times Aut(Q^n)$. By associating a point in $\pi_{P_M}^{-1}(U_\alpha^B)$ with the unique representative, we define the local trivialization of $P_M$ on $U_\alpha^B$ which is denoted by $\phi_\alpha^P : \pi_{P_M}^{-1}(U_\alpha^B) \to U_\alpha^B \times Aut(Q^n)$. Note that the following equality holds for $(b, h) \in U_{\alpha\beta}^B \times Aut(Q^n)$

$$\phi_{\alpha\beta}^P(b, h) := \phi_\alpha^P \circ (\phi_\beta^P)^{-1}(b, h) = (b, \rho_{\alpha\beta} \circ h). \tag{3.2}$$

We denote by $\pi_{\Lambda_M} : \Lambda_M \to B_M$ the associated $\Lambda$-bundle associated with $P_M$.

For each coordinate neighborhood $(U_\alpha^B, \phi_\alpha^B)$ of $B_M$ with $U_\alpha^B \cap \mathcal{S}^{(n-1)} B_M \neq \emptyset$, $Q^n$ acts on the preimage $\pi_{\mathbb{H}^n}^{-1}(\phi_\alpha^B(U_\alpha^B \cap \mathcal{S}^{(n-1)} B_M))$ as the restriction of the regular representation of $Q^n$. For simplicity, we assume that the intersection $U_\alpha^B \cap \mathcal{S}^{(n-1)} B_M$ is connected (otherwise, we make a componentwise argument). Then, $\pi_{\mathbb{H}^n}^{-1}(\phi_\alpha^B(U_\alpha^B \cap \mathcal{S}^{(n-1)} B_M))$ is fixed by a 3-sphere $S_\alpha^3$. We denote by $\mathcal{L}_\alpha$ the rank-one sublattice of $\Lambda$ spanned by the integral element that generates $S_\alpha^3$. Suppose that $(U_\alpha^B, \phi_\alpha^B)$ and $(U_\beta^B, \phi_\beta^B)$ are coordinate neighborhoods satisfying the above conditions and that the intersection $U_{\alpha\beta}^B \cap \mathcal{S}^{(n-1)} B_M$ is nonempty.

By the construction of $\pi_{\Lambda_M} : \Lambda_M \to B_M$, there exists a local trivialization $\varphi_\alpha^B : \pi_{\Lambda_M}^{-1}(U_\alpha^B) \to U_\alpha^B \times \Lambda$ of $\pi_{\Lambda_M} : \Lambda_M \to B_M$ on each $U_\alpha^B$ such that on an overlap $U_{\alpha\beta}^B$ the transition function with respect to $\varphi_\alpha^B$ and $\varphi_\beta^B$ is $\rho_{\alpha\beta}$. We take a subsystem $\{(U_{\alpha_i}^B, \varphi_{\alpha_i}^B)\}_{i \in I}$ of $\{(U_\alpha^B, \varphi_\alpha^B)\}_{\alpha \in A}$ which covers $\mathcal{S}^{(n-1)} B_M$. Notice that under the identification $Aut(Q^n) \cong GL(\Lambda_Q)$, the automorphism $\rho_{\alpha\beta}$ of $Q^n$ in the Definition 2.3 sends $\mathcal{L}_\beta$ isomorphically to $\mathcal{L}_\alpha$. This means that we can obtain the rank-one subbundle $\pi_{\mathcal{L}_M} : \mathcal{L}_M \to \mathcal{S}^{(n-1)} B_M$ of $\pi_{\Lambda_M}|_{\mathcal{S}^{(n-1)} B_M} : \Lambda_M|_{\mathcal{S}^{(n-1)} B_M} \to \mathcal{S}^{(n-1)} B_M$ by setting

$$\mathcal{L}_M := \left(\bigsqcup_i U_{\alpha_i}^B \cap \mathcal{S}^{(n-1)} B_M \times \mathcal{L}_{\alpha_i}\right) / \sim_L, \tag{3.3}$$

where $(b_i, l_i) \in U_{\alpha_i}^B \cap \mathcal{S}^{(n-1)} B_M \times \mathcal{L}_{\alpha_i} \sim_L (b_j, l_j) \in U_{\alpha_j}^B \cap \mathcal{S}^{(n-1)} B_M \times \mathcal{L}_{\alpha_j}$ if and only if $b_i = b_j$ and $l_i = \rho_{\alpha_i \alpha_j}(l_j)$.

**Definition 3.2.** Let $\pi_\mathcal{L} : \mathcal{L} \to \mathcal{S}^{(n-1)} B$ be a rank-one sublattice bundle of the associated lattice bundle $\pi_\Lambda : \Lambda \to B$ defined over the codimension-one stratum $\mathcal{S}^{(n-1)} B$ of a manifold with corners $B$. The bundle $\mathcal{L}$ is said to be *unimodular* if, for every point $b \in \mathcal{S}^{(k)} B$ (the codimension-$k$ stratum), the fibers of $\mathcal{L}$ over a neighborhood of $b$ generate a rank-$k$ direct summand of the full lattice $\Lambda_b \cong \mathbb{Z}^n$.

By construction, it is easy to see that $\pi_{\mathcal{L}_M} : \mathcal{L}_M \to \mathcal{S}^{(n-1)} B_M$ is unimodular. We summarize the above discussion in the following.

**Proposition 3.3.** *Associated with a local $Q^n$-action $\mathcal{Q}$ on $M^{4n}$, there exists a characteristic pair $(P_M, \mathcal{L}_M)$, where $P_M$ and $\mathcal{L}_M$ are defined by (3.1) and (3.3), respectively.*

**Example 3.4** (Quoric manifolds)**.** A quoric manifold $M^{4n}$ admits a globally defined regular $Q^n$-action whose orbit space is a simple polytope $P^n$ (see [Hopkinson, 2012]). The associated characteristic pair $(P^n, \lambda)$ is by definition unimodular and the Čech cocycle is trivial. In other words, the local action globalizes.

**Definition 3.5.** An *isomorphism* $f_P : (P_1, \mathcal{L}_1) \to (P_2, \mathcal{L}_2)$ between characteristic pairs is a bundle isomorphism $f_P : P_1 \to P_2$ which covers a stratification preserving homeomorphism $f_B : B_1 \to B_2$ such that the induced bundle isomorphism $f_\Lambda : \Lambda_1 \to \Lambda_2$ induced by $f_P$ sends $\mathcal{L}_1$ isomorphically to $\mathcal{L}_2$. $(P_1, \mathcal{L}_1)$ and $(P_2, \mathcal{L}_2)$ are *isomorphic* if there exists an isomorphism between them.

The isomorphism class of the characteristic pair $(P_M, \mathcal{L}_M)$ is a combinatorial invariant of a local $Q^n$-action on $M^{4n}$.

**Lemma 3.6.** *For $i = 1, 2$, let $(M_i^{4n}, \mathcal{Q}_i)$ be a $4n$-dimensional manifold with a local $Q^n$-action. If there is a homeomorphism $f_M : (M_1^{4n}, \mathcal{Q}_1) \to (M_2^{4n}, \mathcal{Q}_2)$, then $f_M$ induces an isomorphism $f_{P_M} : (P_{M_1^{4n}}, \mathcal{L}_{M_1^{4n}}) \to (P_{M_2^{4n}}, \mathcal{L}_{M_2^{4n}})$ between characteristic pairs associated with $M_1^{4n}$ and $M_2^{4n}$.*

*Proof.* For simplicity, we denote $M_1^{4n}$ and $M_2^{4n}$ by $M_1$ and $M_2$, respectively. Let now $\{(U_\beta^{M_1}, \varphi_\beta^{M_1})\}_{\beta \in \mathcal{B}} \in \mathcal{Q}_1$ and $\{(U_\alpha^{M_2}, \varphi_\alpha^{M_2})\}_{\alpha \in \mathcal{A}} \in \mathcal{Q}_2$ be maximal weakly regular atlases of $M_1$ and $M_2$, and $\{(U_\beta^{B_1}, \phi_\beta^{B_1})\}_{\beta \in \mathcal{B}}$ and $\{(U_\alpha^{B_2}, \phi_\alpha^{B_2})\}_{\alpha \in \mathcal{A}}$ be atlases of $B_{M_1}$ and $B_{M_2}$ induced by $\{(U_\beta^{M_1}, \varphi_\beta^{M_1})\}_{\beta \in \mathcal{B}}$ and $\{(U_\alpha^{M_2}, \varphi_\alpha^{M_2})\}_{\alpha \in \mathcal{A}}$, respectively. Suppose that $f_M : (M_1, \mathcal{Q}_1) \to (M_2, \mathcal{Q}_2)$ is a homeomorphism and $f_B : B_{M_1} \to B_{M_2}$ is the homeomorphism induced by $f_M$. By definition, on each nonempty overlap $U_\beta^{B_1} \cap f_B^{-1}(U_\alpha^{B_2})$, there exists an automorphism $\rho_{\alpha\beta}^f$ of $Q^n$ such that $\varphi_\alpha^{M_2} \circ f_M \circ (\varphi_\beta^{M_1})^{-1}$ is $\rho_{\alpha\beta}^f$-equivariant. It is easy to see that the equality

$$\rho_{\alpha_0\beta_0}^f \circ \rho_{\beta_0\beta_1}^{M_1} = \rho_{\alpha_0\alpha_1}^{M_2} \circ \rho_{\alpha_1\beta_1}^f \tag{3.4}$$

holds on a nonempty intersection $U_{\beta_0\beta_1}^{B_1} \cap f_B^{-1}(U_{\alpha_0\alpha_1}^{B_2})$, where $\rho_{\beta_0\beta_1}^{M_1}$ and $\rho_{\alpha_0\alpha_1}^{M_2}$ are the automorphisms of $Q^n$ in Definition 2.3 with respect to $M_1$ and $M_2$, respectively. We define a bundle isomorphism $(f_P)_{\alpha\beta} : U_\beta^{B_1} \cap f_B^{-1}(U_\alpha^{B_2}) \times \text{Aut}(Q^n) \to U_\alpha^{B_2} \cap f_B^{-1}(U_\beta^{B_1}) \times \text{Aut}(Q^n)$ by

$$(f_P)_{\alpha\beta}(b, h) := (f_B(b), \rho_{\alpha\beta}^f \circ h).$$

Then, (3.4) implies that the equation

$$(f_P)_{\alpha_0\beta_0} \circ \varphi_{\beta_0\beta_1}^{P_{M_1}} = \varphi_{\alpha_0\alpha_1}^{P_{M_2}} \circ (f_P)_{\alpha_1\beta_1}$$

holds on $U_{\beta_0\beta_1}^{B_1} \cap f_B^{-1}(U_{\alpha_0\alpha_1}^{B_2})$, where $\varphi_{\beta_0\beta_1}^{P_{M_1}}$ (resp. $\varphi_{\alpha_0\alpha_1}^{P_{M_2}}$) is the overlap map defined by (3.2) for $\pi_{P_{M_1}} : P_{M_1} \to B_1$ (resp. $\pi_{P_{M_2}} : P_{M_2} \to B_2$) on $U_{\beta_0\beta_1}^{B_1}$ (resp. $U_{\alpha_0\alpha_1}^{B_2}$). Therefore, we can patch them together to obtain the bundle isomorphism $f_P : P_{M_1} \to P_{M_2}$ which covers $f_B$. □

### 3.3 Canonical model

In this section we construct the canonical model by using combinatorial data, specifically a characteristic pair $(P, \mathcal{L})$, defined on a base space $B$, an n-manifold with corners. Denote again by $\pi_Q : Q_P \to B$ the quaternionic torus bundle associated to $P$ by the natural action of $\text{Aut}(Q^n)$ on $Q^n$, whose fibers are $Q^n$. First we shall explain that for any $k$-dimensional part $\mathcal{S}^{(k)}B$, $(P, \mathcal{L})$ determines a rank-$(n-k)$ subtorus bundle of the restriction of $\pi_Q : Q_P \to B$ to $\mathcal{S}^{(k)}B$. Let $\{U_\alpha^B\}$ be an open covering of $B$ such that on each $U_\alpha^B$ there exists a local trivialization $\varphi_\alpha^P : \pi_Q^{-1}(U_\alpha^B) \to U_\alpha^B \times \mathrm{Aut}(Q^n)$. On each nonempty overlap $U_{\alpha\beta}^B$ we denote by $\rho_{\alpha\beta}$ the transition function with respect to $\varphi_\alpha^P$ and $\varphi_\beta^P$, namely,

$$\varphi_\alpha^P \circ (\varphi_\beta^P)^{-1}(b, f) = (b, \rho_{\alpha\beta} f)$$

for $(b, f) \in U_\beta^B \times \mathrm{Aut}(Q^n)$. Note that $\rho_{\alpha\beta}$ is locally constant since $\mathrm{Aut}(Q^n)$ is discrete. For simplicity we assume that each $U_{\alpha\beta}^B$ is connected so that $\rho_{\alpha\beta}$ can be thought of as an element of $\mathrm{Aut}(Q^n)$. The map $\varphi_\alpha^P$ induces local trivializations of the associated bundles $Q_P$ and $\Lambda_P$ denoted by $\varphi_\alpha^Q : \pi_Q^{-1}(U_\alpha^B) \to U_\alpha^B \times Q^n$ and $\varphi_\alpha^\Lambda : \pi_\Lambda^{-1}(U_\alpha^B) \to U_\alpha^B \times \Lambda$, respectively.

For $\mathcal{S}^{(k)}B$ we take $U_\alpha^B$ with $U_\alpha^B \cap \mathcal{S}^{(k)}B \neq \emptyset$. Replacing $U_\alpha^B$ by a sufficiently small one if necessary, we may assume that the intersection $U_\alpha^B \cap \mathcal{S}^{(n-1)}B$ has exactly $n - k$ connected components, say $(U_\alpha^B \cap \mathcal{S}^{(n-1)}B)_1, \ldots, (U_\alpha^B \cap \mathcal{S}^{(n-1)}B)_{n-k}$. For $k = n$, this means that $U_\alpha^B$ is contained in $\mathcal{S}^{(n)}B$.

For $k < n$, there are $n - k$ rank-one sublattices $L_1, \ldots, L_{n-k}$ of $\Lambda$ such that for $a = 1, \ldots, n-k$, $\varphi_\alpha^\Lambda$ sends the restriction of $\pi_\Lambda : \mathcal{L} \to \mathcal{S}^{(n-1)}B$ to $(U_\alpha^B \cap \mathcal{S}^{(n-1)}B)_a$ isomorphically to the trivial rank-one subbundle $(U_\alpha^B \cap \mathcal{S}^{(n-1)}B)_a \times L_a$ of $(U_\alpha^B \cap \mathcal{S}^{(n-1)}B)_a \times \Lambda$. Since $\mathcal{L}$ is unimodular, $L_1, \ldots, L_{n-k}$ generate the isotropy subgroups at lower-dimensional strata (e.g., corners or boundaries of $B$) which are quaternionic tori of dimension $n - k$. That is,

$$Q_{U_\alpha^B \cap \mathcal{S}^{(k)}B} \cong (S^3)^{n-k} = Q^{n-k}.$$

For $k = n$, we define $Q_{U_\alpha^B \cap \mathcal{S}^{(n)}B}$ to be the trivial subgroup. Note that when $(P, \mathcal{L})$, $\{U_\alpha^B\}$ and $\varphi_\alpha^P$ are induced by some local $Q^n$-action $\mathcal{Q}$ on $M^{4n}$, $Q_{U_\alpha^B \cap \mathcal{S}^{(k)}B_M}$ is the common $(n-k)$-dimensional stabilizer (isotropy subgroup) of the $Q^n$-action on $M^{4n}$.

Suppose that another local trivialization $\{U_\beta^B\}$ satisfies $U_\beta^B \cap \mathcal{S}^{(k)}B \neq \emptyset$. By definition of $(P, \mathcal{L})$, the automorphism $\rho_{\alpha\beta}$ sends $Q_{U_\alpha^B \cap \mathcal{S}^{(k)}B}$ isomorphically to $Q_{U_\beta^B \cap \mathcal{S}^{(k)}B}$. Hence, patching them together using $\rho_{\alpha\beta}$, we obtain a subtorus bundle, denoted by

$$\pi_{Q_{\mathcal{S}^{(k)}B}} : Q_{\mathcal{S}^{(k)}B} \to \mathcal{S}^{(k)}B.$$

**Definition 3.7.** Two elements $q, q' \in Q_P$ are said to be *equivalent*, writing $q \sim q'$, if $\pi_Q(q) = \pi_Q(q')$ and $q'q^{-1} \in \pi_{Z_{\mathcal{S}^{(k)}B}}^{-1}(\pi_Q(q))$, provided that $\pi_Q(q) \in \mathcal{S}^{(k)}B$.

We denote by $M_{(P,\mathcal{L})}$ the quotient space of $Q_P$ by this equivalence relation, with a natural local $Q^n$-action, whose orbit space is $B$ and whose orbit map is $\pi_{M_{(P,\mathcal{L})}}$.

**Definition 3.8.** We call $M_{(P,\mathcal{L})}$ the *canonical model* associated with $(P, \mathcal{L})$. In particular, when $(P, \mathcal{L})$ is the characteristic pair $(P_M, \mathcal{L}_M)$ of a local $Q^n$-action $\mathcal{Q}$ on a $4n$-dimensional manifold $M^{4n}$, we also call $M_{(P_M, \mathcal{L}_M)}$ the *canonical model associated to* $(M^{4n}, \mathcal{Q})$.

We give some properties of a canonical model.

**Lemma 3.9.** *For a characteristic pair $(P, \mathcal{L})$, the projection $\pi_{M_{(P,\mathcal{L})}} : M_{(P,\mathcal{L})} \to B$ admits a section*

$$s_{M_{(P,\mathcal{L})}} : B \to M_{(P,\mathcal{L})}.$$

*Proof.* By Definition 3.7, $M_{(P,\mathcal{L})}$ is the quotient of the torus bundle $Q_P \to B$ by the equivalence relation that identifies fibers over $\mathcal{S}^{(k)}B \subset B$ along the subgroups $Q_{\mathcal{S}^{(k)}B} \subset Q^n$. Each fiber of $Q_P \to B$, and hence of $M_{(P,\mathcal{L})} \to B$, is modeled on the group $Q^n$, which is a topological group.

Since $Q_P \to B$ is locally trivial and its structure group is the discrete group $\mathrm{Aut}(Q^n)$, the bundle admits local trivializations over an open cover $\{U_\alpha^B\}$ of $B$. Within each local trivialization $\varphi_\alpha^P : \pi_P^{-1}(U_\alpha^B) \to U_\alpha^B \times Q^n$, we may define a local section by mapping each $b \in U_\alpha^B$ to the identity element $e \in Q^n$, i.e., $s'_\alpha(b) = \varphi_\alpha^{P-1}(b, e) \in Q_P$.

These local sections are compatible under the quotient that defines $M_{(P,\mathcal{L})}$, since the identity element of $Q^n$ lies in every isotropy subgroup $Q_{\mathcal{S}^{(k)}B}$ and is fixed under the action of the transition functions $\rho_{\alpha\beta} \in \mathrm{Aut}(Q^n)$. Therefore, the composition of $s' : B \to Q_P$ with the natural projection $Q_P \to M_{(P,\mathcal{L})}$ yields a section

$$s_{M_{(P,\mathcal{L})}} : B \to M_{(P,\mathcal{L})}$$

satisfying $\pi_{M_{(P,\mathcal{L})}} \circ s_{M_{(P,\mathcal{L})}} = \mathrm{id}_B$. □

The next result follows directly from the construction of the local trivializations used in defining the canonical model.

**Lemma 3.10.** *For $i = 1, 2$, let $B_i$ be an $n$-dimensional manifold with corners and $(P_i, L_i)$ a characteristic pair on $B_i$. Then, an isomorphism $f_P : (P_1, L_1) \to (P_2, L_2)$ induces a homeomorphism $f_{M_{(P,L)}} : M_{(P_1,L_1)} \to M_{(P_2,L_2)}$ between the canonical models of $(P_1, L_1)$ and $(P_2, L_2)$.*

*Proof.* The canonical model $M_{(P,\mathcal{L})}$ is constructed from the associated torus bundle $Q_P \to B$, together with the isotropy data given by $\mathcal{L}$. An isomorphism $f_P : (P_1, \mathcal{L}_1) \to (P_2, \mathcal{L}_2)$ consists of a homeomorphism $f_B : B_1 \to B_2$ and a bundle isomorphism $Q_{P_1} \to Q_{P_2}$ covering $f_B$, which also respects the sublattice data encoded by $\mathcal{L}_i$.

This isomorphism preserves the local trivializations of the torus bundles and maps isotropy subtori of $Q^n$ over $\mathcal{S}^{(k)}B_1$ isomorphically to those over $\mathcal{S}^{(k)}B_2$. Consequently, the equivalence relations used in the construction of $M_{(P_1,\mathcal{L}_1)}$ and $M_{(P_2,\mathcal{L}_2)}$ are compatible under this isomorphism. Therefore, the map descends to a homeomorphism of canonical models

$$f_{M_{(P,\mathcal{L})}} : M_{(P_1,\mathcal{L}_1)} \to M_{(P_2,\mathcal{L}_2)}.$$

□

## 4 Cohomological obstruction and topological classification

Let $(M^{4n}, \mathcal{Q})$ be a $4n$-dimensional manifold equipped with a local $Q^n$-action. In this section, we investigate the role of cohomological invariants in the classification of such manifolds. As a first step, we compare $M^{4n}$ with its associated canonical model $M_{(P,\mathcal{L})}$. This is achieved by a cohomological invariant known as the Euler class. In particular, to the orbit map we associate its Euler class and show that its vanishing is equivalent to the existence of a section of the orbit map (see [Hatcher, 2017, Section 3.3]). Then, using this invariant, we topologically classify $4n$-dimensional manifolds equipped with a local $Q^n$-action (local quaternionic torus actions). We follow closely [Yoshida, 2011] and [Gkeneralis and Prassidis, 2025].

### 4.1 The Euler Class of a Local $Q^n$-Action

Notice that, by construction, the orbit map

$$\pi_{\mathbb{H}^n} : \mathbb{H}^n \to \mathbb{R}^n_{\geq}, \ \ \pi_{\mathbb{H}^n}(h) := (|h_1|^4, ..., |h_n|^4)$$

admits a section $i : \mathbb{R}^n_{\geq} \to \mathbb{H}^n$ defined by $i(h) = (h_1^{1/4}, \cdots, h_n^{1/4})$.

**Proposition 4.1.** *Let $(M^{4n}, \mathcal{Q})$ be a $4n$-dimensional manifold equipped with a local $Q^n$-action. If $\pi_M : M^{4n} \to B_M$ has a section $s$ and $i$ is the section of $\pi_{\mathbb{H}^n}$, then there exists a weakly regular atlas of $M$, $\{(U_\alpha^M, \phi_\alpha^M)\}_{\alpha \in \mathcal{A}}$ such that the following diagram commutes*

$$\begin{array}{ccc} U_\alpha^M & \xrightarrow{\phi_\alpha^M} & \phi_\alpha^M(U_\alpha^M) \\ s\uparrow & & i\uparrow \\ U_\alpha^B & \xrightarrow{\phi_\alpha^B} & \phi_\alpha^B(U_\alpha^B) \end{array}$$

where $\{(U_\alpha^B, \phi_\alpha^B)\}_{\alpha \in \mathcal{A}}$ is the atlas of $B_M$ induced by $\{(U_\alpha^M, \phi_\alpha^M)\}_{\alpha \in \mathcal{A}}$.

*Proof.* Let $\{(U_\alpha^M, \psi_\alpha^M)\}_{\alpha \in \mathcal{A}} \in \mathcal{Q}$ be a weakly regular atlas of $M^{4n}$. First, we construct an atlas of $B_M$ as follows. We put $U_\alpha^B := U_\alpha^M/Q^n$. Now, we restrict $\pi_{\mathbb{H}^n}$ to $\psi_\alpha^M(U_\alpha^M)$, which induces a homeomorphism from $U_\alpha^B$ to the open subset $\pi_{\mathbb{H}^n}(\psi_\alpha^M(U_\alpha^M))$ of $\mathbb{R}^n_{\geq}$, which is denoted by $\psi_\alpha^B(U_\alpha^B)$. Thus $\{(U_\alpha^B, \psi_\alpha^B)\}_{\alpha \in \mathcal{A}}$ is the atlas of $B_M$ induced by $\{(U_\alpha^M, \psi_\alpha^M)\}_{\alpha \in \mathcal{A}}$. By construction, this means that the following diagram commutes

$$\begin{array}{ccc} U_\alpha^M \subset M^{4n} & \xrightarrow{\psi_\alpha^M} & \psi_\alpha^M(U_\alpha^M) \subset \mathbb{H}^n \\ \downarrow{\pi_M} & & \downarrow{\pi_{\mathbb{H}^n}} \\ U_\alpha^B \subset B_M & \xrightarrow{\psi_\alpha^B} & \psi_\alpha^B(U_\alpha^B) \subset \mathbb{R}^n_{\geq} \end{array}$$

For each $\alpha \in \mathcal{A}$ and $b \in U_\alpha^B$, we have that $\pi_M \circ (\psi_a^M)^{-1} = (\psi_a^B)^{-1} \circ \pi_{\mathbb{H}^n}$. Then

$$\pi_M \circ (\psi_a^M)^{-1} \circ i = (\psi_a^B)^{-1} \circ \pi_{\mathbb{H}^n} \circ i = (\psi_a^B)^{-1}.$$

Thus for $b \in U_\alpha^B$, $\pi_M \circ (\psi_a^M)^{-1} \circ i \circ \psi_a^B(b) = b = \pi_M \circ s(b)$. The equality

$$\theta_\alpha(b) \cdot s(b) = (\psi_\alpha^M)^{-1} \circ i \circ \psi_\alpha^B(b)$$

for $b \in U_\alpha^B$ determines a local section in the canonical model $M_{(P,\mathcal{L})}$, $\theta_\alpha : U_\alpha^B \to Q^n$.

As in [Yoshida, 2011, Proposition 5.1], we define a new coordinate system $\phi_\alpha^M$ on $U_\alpha^M$: for $x \in \pi_M^{-1}(U_\alpha^B)$, set

$$\phi_\alpha^M(x) := \psi_\alpha^M(\theta_\alpha(\pi_M(x)) \cdot x),$$

which satisfies $\phi_\alpha^M(s(b)) = \psi_\alpha^M(\theta_\alpha(b) \cdot s(b)) = \psi_\alpha^M((\psi_\alpha^M)^{-1} \circ i \circ \psi_\alpha^B(b)) = i \circ \psi_\alpha^B(b)$. In other words, the following diagram commutes

$$\begin{array}{ccc} U_\alpha^M & \xrightarrow{\phi_\alpha^M} & \phi_\alpha^M(U_\alpha^M) \\ s\uparrow & & i\uparrow \\ U_\alpha^B & \xrightarrow{\phi_\alpha^B} & \phi_\alpha^B(U_\alpha^B) \end{array}$$

Therefore, the collection $\{(U_\alpha^M, \phi_\alpha^M)\}_{\alpha \in \mathcal{A}}$ is the desired weakly regular atlas on $M^{4n}$, and the associated base charts remain unchanged, i.e., $\phi_\alpha^B := \psi_\alpha^B$. □

**Proposition 4.2.** *The projection $\pi_M : M^{4n} \to B_M$ admits a section if and only if there exist a homeomorphism between $(M^{4n}, \mathcal{Q})$ and $M_{(P_M, \mathcal{L}_M)}$ which covers the identity.*

*Proof.* Let $h$ denote the $Q^n$-homeomorphism between $M^{4n}$ and $M_{(P_M, \mathcal{L}_M)}$. Define $s : B_M \to M^{4n}$ such that,

$$s(x) = h^{-1} \circ s_{M_{(P,\mathcal{L})}}(x),$$

where $s_{M_{(P,\mathcal{L})}}$ is defined in Lemma 3.9. Then $\pi \circ s(x) = \pi \circ h^{-1} \circ s_{M_{(P,\mathcal{L})}}(x) = \pi_{M_{(P,\mathcal{L})}} \circ s_{M_{(P,\mathcal{L})}}(x) = id_P(x) = x$.

On the other hand, let $s : B_M \to M^{4n}$ be a section of $\pi_M$. Define the map $h : M_{(P_M, \mathcal{L}_M)} \to M^{4n}$ by $h(q, x) = q \cdot s(x)$ which is obviously a $Q^n$-homeomorphism. □

As in classical bundle theory, we now define a cohomological invariant of a local $Q^n$-action, known as the Euler class.

Let $(M^{4n}, \mathcal{Q})$ be a $4n$-dimensional manifold equipped with a local $Q^n$-action and associated orbit map $\pi_M : M^{4n} \to B_M$. Let $\{(U_\alpha^M, \varphi_\alpha^M)\}_{\alpha \in \mathcal{A}} \in \mathcal{Q}$ be a weakly regular atlas of $M^{4n}$. Let $U_\alpha^B := \pi_M(U_\alpha^M)$ and $\{U_\alpha^B\}_{\alpha \in \mathcal{A}}$ be the induced open cover of $B_M$.

Suppose that $\pi_M$ admits a section $s : B_M \to M^{4n}$, and let $h_\alpha : U_\alpha^B \times Q^n/\sim \,\to\, U_\alpha^M$ be the local trivialization defined by

$$h_\alpha(q, b) := q \cdot s(b).$$

On overlaps $U_{\alpha\beta}^B := U_\alpha^B \cap U_\beta^B$, both $h_\alpha$ and $h_\beta$ trivialize $\pi_M$, and any point $x \in \pi_M^{-1}(U_{\alpha\beta}^B)$ satisfies

$$h_\alpha \circ h_\beta^{-1}(x) = \theta_{\alpha\beta}^M(b) \cdot x$$

for some uniquely determined $\theta_{\alpha\beta}^M(b) \in Q^n$.

The next result ensures that this description is consistent.

**Proposition 4.3.** *For each $b \in B_M$, the action of $Q^n$ on the fiber $\pi_M^{-1}(b)$ is simply transitive. Hence, for any pair of local trivializations $h_\alpha, h_\beta$ over $U_{\alpha\beta}^B$, there exists a unique group element $\theta_{\alpha\beta}^M(b) \in Q^n$ such that*

$$h_\alpha \circ h_\beta^{-1}(x) = \theta_{\alpha\beta}^M(b) \cdot x.$$

**Remark 4.4.** The above Proposition 4.3 is the quaternionic analogue of [Yoshida, 2011, Prop. 4.12]. The proof is similar in spirit, relying on the fact that the local $Q^n$-actions on the fibers are free and transitive and thus is omitted.

Thus, the family $\check{\theta} = \{\theta_{\alpha\beta}^M\}_{\alpha,\beta}$ defines a Čech 1-cochain with values in the sheaf $\mathscr{S}$ of local sections of the bundle $\pi_M : M^{4n} \to B_M$.

The cocycle conditions follow immediately from associativity and the definition of the $\theta_{\alpha\beta}^M$:

$$\theta_{\alpha\alpha}^M = \text{id},$$
$$\theta_{\alpha\beta}^M \cdot \theta_{\beta\gamma}^M \cdot \theta_{\gamma\alpha}^M = \text{id} \quad \text{on } U_{\alpha\beta\gamma}^B.$$

Hence, $\check{\theta}$ defines a Čech 1-cocycle with values in $\mathscr{S}$, and determines a cohomology class in the Čech cohomology group $\check{H}^1(B_M, \mathscr{S})$.

**Definition 4.5.** The cohomology class $e(M) \in \check{H}^1(B_M, \mathscr{S})$ defined by the cocycle $\check{\theta} = \{\theta_{\alpha\beta}^M\}$ is called the *Euler class* of the bundle $\pi_M : M^{4n} \to B_M$.

**Theorem.** *The orbit map $\pi : M^{4n} \to B_M$ has a global section if and only if $e(M)$ vanishes.*

*Proof.* Let $\delta$ denote the Čech differential. The class $e(M)$ vanishes when a representative $\check{\theta} \in e(M)$ is cohomologous to the zero cocycle. Equivalently, there exists a 0 - cocycle $\theta$, such that $\delta\theta = \check{\theta}$. But, here, cocycles $\theta \in \check{H}^0(B_M, \mathscr{S})$ are cochains $\{\theta_\alpha\}$ which satisfy

$$\theta_\alpha = \check{\theta}_{\alpha\beta}\theta_\beta \tag{4.1}$$

for all $x \in U_{\alpha\beta}$. The equation (4.1) becomes

$$\theta_\alpha(b)x = \theta_{\alpha\beta}^M \theta_\beta(b)x = h_\alpha \circ h_\beta^{-1}(x)\theta_\beta(b)x \Leftrightarrow h_\alpha^{-1}(\theta_\alpha(b)x) = h_\beta^{-1}(\theta_\beta(b)x)$$

for any $b \in U_{\alpha\beta}^B$, where we used the fact that $h_\alpha$ and $h_\beta$ are equivariant maps. In other words, $h_\alpha$ and $h_\beta$ agree on each non empty overlap $U_{\alpha\beta}$. This means that the map $h : M^{4n} \to M_{(P,\mathcal{L})}$ defined by $h(x) = h_\alpha^{-1}(\theta_\alpha(b)x)$ for $x \in U_\alpha^M$ and $\pi(x) = b \in U_\alpha^B$ is a $Q^n$-homeomorphism if and only if $e(M)$ vanishes.

So if the orbit map $\pi : M^{4n} \to B_M$ has a section, by Proposition 4.2, $M^{4n}$ and $M_{(P,\mathcal{L})}$ are $Q^n$-homeomorphic and thus $e(M) = 0$.

On the other hand, if the Euler class $e(M)$ vanishes then we have that $M^{4n}$ and $M_{(P,\mathcal{L})}$ are $Q^n$-homeomorphic. Again by Proposition 4.2 the orbit map $\pi : M^{4n} \to B_M$ has a section. □

### 4.2 The topological classification

In this subsection, we establish a topological classification of compact $4n$-dimensional manifolds equipped with local $Q^n$-actions. We show that such manifolds are completely determined, up to equivariant homeomorphism, by their characteristic pair and the Euler class.

**Theorem.** *Let $(M_1^{4n}, \mathcal{Q}_1)$ and $(M_2^{4n}, \mathcal{Q}_2)$ be two compact $4n$-dimensional manifolds equipped with local $Q^n$-actions. Then $(M_1, \mathcal{Q}_1)$ and $(M_2, \mathcal{Q}_2)$ are equivariantly homeomorphic if and only if there exists an isomorphism of characteristic pairs $f_{P_M} : (P_{M_1}, \mathcal{L}_{M_1}) \xrightarrow{\cong} (P_{M_2}, \mathcal{L}_{M_2})$ such that $f_{P_M}^* e(M_2) = e(M_1)$.*

*Moreover, for any characteristic pair $(P, \mathcal{L})$ on a compact $n$-dimensional manifold with corners $B$, and any cohomology class $e \in \check{H}^1(B, \mathcal{S})$, there exists a compact $4n$-dimensional manifold $M^{4n}$ equipped with a local $Q^n$-action $\mathcal{Q}$, such that its characteristic pair is $(P, \mathcal{L})$ and its Euler class is $e$.*

*Proof.* Suppose that there exists an equivariant homeomorphism $f_M : (M_1^{4n}, \mathcal{Q}_1) \to (M_2^{4n}, \mathcal{Q}_2)$. By Lemma 3.6, $f_M$ induces a bundle isomorphism

$$f_{P_M} : (P_{M_1}, \mathcal{L}_{M_1}) \to (P_{M_2}, \mathcal{L}_{M_2})$$

between the characteristic pairs. Moreover, since $f_M$ is $Q^n$-equivariant and respects the local trivializations defining the Euler cocycles, it pulls back the cocycle $\check{\theta}_{M_2}$ representing $e(M_2)$ to a cocycle cohomologous to $\check{\theta}_{M_1}$. Therefore, we have

$$f_{P_M}^* e(M_2) = e(M_1)$$

in $\check{H}^1(B_{M_1}, \mathcal{S})$, as desired.

On the other hand, suppose that there exists an isomorphism of characteristic pairs $f_{P_M} : (P_{M_1}, \mathcal{L}_{M_1}) \to (P_{M_2}, \mathcal{L}_{M_2})$ such that $f_{P_M}^* e(M_2) = e(M_1)$. Let $f_B : B_{M_1} \to B_{M_2}$ be the induced stratification-preserving homeomorphism covered by $f_{P_M}$.

Let $\{(U_\beta^1, \varphi_\beta^1)\} \in \mathcal{Q}_1$, $\{(U_\alpha^2, \varphi_\alpha^2)\} \in \mathcal{Q}_2$ be weakly regular atlases for $M_1$ and $M_2$, respectively, and let $\{(U_\beta^{B_1}, \varphi_\beta^{B_1})\}$, $\{(U_\alpha^{B_2}, \varphi_\alpha^{B_2})\}$ be the induced atlases on $B_{M_1}, B_{M_2}$. For each chart, we can take homeomorphisms

$$h_\beta^1 : \pi_{M_1}^{-1}(U_\beta^{B_1}) \to \pi_{(P_{M_1}, \mathcal{L}_{M_1})}^{-1}(U_\beta^{B_1}), \quad h_\alpha^2 : \pi_{M_2}^{-1}(U_\alpha^{B_2}) \to \pi_{(P_{M_2}, \mathcal{L}_{M_2})}^{-1}(U_\alpha^{B_2}),$$

satisfying the appropriate local trivialization compatibility for $M_1$ and $M_2$, as defined in Proposition 4.2.

By the assumption $f_P^* e(M_2) = e(M_1)$, and refining the covers if needed, we may assume that there exists a Čech 0-cochain $\{\theta_\beta\}$ such that the relation

$$\theta_{\beta_0 \beta_1}^{M_1}(b) = f_{P_M}^{-1}\left(\theta_{\alpha_0 \alpha_1}^{M_2}(f_B(b))\right) \cdot \theta_{\beta_1}(b) \cdot \theta_{\beta_0}(b)^{-1} \tag{4.2}$$

holds on each non-empty intersection $b \in U_{\beta_0 \beta_1}^{B_1} \cap f_B^{-1}(U_{\alpha_0 \alpha_1}^{B_2})$.

Now define, for each such $b$, the local homeomorphism

$$f_{\alpha\beta} : \pi_{M_1}^{-1}\left(U_\beta^{B_1} \cap f_B^{-1}(U_\alpha^{B_2})\right) \to \pi_{M_2}^{-1}\left(f_B(U_\beta^{B_1}) \cap U_\alpha^{B_2}\right)$$

by
$$f_{\alpha\beta}(x) := (h_\alpha^2)^{-1} \circ f_{P_M}\left(\theta_\beta(\pi_{M_1}(x)) \cdot h_\beta^1(x)\right).$$

This construction gives local homeomorphisms that patch together on overlaps using Equation (4.2). Therefore, they define a global $Q^n$-equivariant homeomorphism

$$f_M : (M_1, \mathcal{Q}_1) \to (M_2, \mathcal{Q}_2),$$

covering $f_B$, as required.

Suppose that $(P, \mathcal{L})$ is a characteristic pair on an $n$-dimensional manifold $B$ with corners, and let $e \in \check{H}^1(B; \mathcal{S})$ be a cohomology class. Take a representative $\{\theta_{\alpha\beta}\}$ of $e$ with respect to a sufficiently fine open cover $\{U_\alpha\}$ of $B$.

We construct a new $4n$-dimensional manifold $(M, \mathcal{Q})$ equipped with a local $Q^n$-action by setting

$$M := \left(\bigsqcup_\alpha \pi_{M_{(P,\mathcal{L})}}^{-1}(U_\alpha^B)\right) \Big/ \sim,$$

where for $x_\alpha \in \pi_{M_{(P,\mathcal{L})}}^{-1}(U_\alpha^B)$ and $x_\beta \in \pi_{M_{(P,\mathcal{L})}}^{-1}(U_\beta^B)$, we define $x_\alpha \sim x_\beta$ if and only if $\pi_{M_{(P,\mathcal{L})}}(x_\alpha) = \pi_{M_{(P,\mathcal{L})}}(x_\beta)$ and $x_\alpha = \theta_{\alpha\beta}(\pi_{M_{(P,\mathcal{L})}}(x_\alpha)) \cdot x_\beta$.

It is straightforward to check that the resulting manifold $M$ carries a local $Q^n$-action whose characteristic pair and Euler class are equal to $(P, \mathcal{L})$ and $e$, respectively. □

**Corollary 4.6.** *Let $(M_1^{4n}, \mathcal{Q}_1)$ and $(M_2^{4n}, \mathcal{Q}_2)$ be compact $4n$-dimensional manifolds with local $Q^n$-actions. If their characteristic pairs and Euler classes are isomorphic via*

$$f_P : (P_{M_1}, \mathcal{L}_{M_1}) \xrightarrow{\cong} (P_{M_2}, \mathcal{L}_{M_2}), \quad \text{with} \quad f_P^* e(M_2) = e(M_1),$$

*then any two such local $Q^n$-manifolds are homeomorphic as local $Q^n$-spaces.*

*In particular, $(P_M, \mathcal{L}_M, e(M))$ completely determine the topological type of $(M^{4n}, \mathcal{Q})$ up to homeomorphism.*

The above theorem can also be applied in the case of $4n$-dimensional manifolds equipped with a locally regular $Q^n$-action (as in [Yoshida, 2011, Theorem 6.2]).

**Corollary 4.7.** *Locally regular quaternionic torus actions are classified by the characteristic pair and the Euler class of the orbit map up to equivariant homeomorphisms.*

**Remark 4.8.** Corollary 4.7 is a generalization of the topological classification of quoric manifolds by Hopkinson [Hopkinson, 2012].

## 5   Tetraplectic geometry

In this section, we investigate the following question:

- *Which local quaternionic torus actions can be equipped with tetraplectic forms so as to admit local quaternionic toric fibrations?*

### 5.1   Basic notions

In this section, we review basic facts on tetraplectic structure, following the work of Foth; for details see [Foth, 2002]. This background will be essential for our study of a particular subclass of manifolds equipped with local quaternionic torus actions that admit tetraplectic structures. In what follows the Lie algebra of $Q = SU(2)$ will be denoted by $\mathfrak{q} = \mathfrak{su}(2)$.

**Definition 5.1.** Let $M^{4n}$ be a real manifold. A 4-form $\psi$ on $M^{4n}$ is said to be *tetraplectic* if

(1) $\psi$ is closed, i.e. $d\psi = 0$;

(2) $\psi$ is non-degenerate, i.e. the map $v \mapsto \iota_v \psi$ that contracts $\psi$ along a tangent vector field $v$ has trivial kernel.

We also call $\psi$ a *tetraplecic structure* and the pair $(M^{4n}, \psi)$ is called a *tetraplectic manifold*.

**Definition 5.2.** Let $(M^{4n}, \psi)$ and $(X^{4n}, \psi')$ be smooth $4n$-dimensional manifolds, each one equipped with a tetraplectic 4-form. A diffeomorphism $f : M^{4n} \to X^{4n}$ is called a *tetraplectomorphism* if it preserves the tetraplectic structures, i.e.,

$$f^*\psi' = \psi.$$

**Example 5.3** (Symplectic manifolds)**.** A natural class of examples of tetraplectic manifolds is given by $4n$-dimensional symplectic manifolds. Indeed starting from a manifold with symplectic form $\omega$ we obtain a tetraplectic manifold by endowing it with the 4-form $\omega \wedge \omega$.

The following example is fundamental in this paper.

**Example 5.4** (The space of quaternions)**.** Let $\mathbb{H}^n \cong \mathbb{R}^{4n}$ be the space of quaternions. The standard tetraplectic structure is

$$\psi_{\mathbb{H}^n} = \sum_{i=1}^n dx_{4i-3} \wedge dx_{4i-2} \wedge dx_{4i-1} \wedge dx_{4i}$$

where $x_1, ..., x_{4n}$ is the coordinate system on $\mathbb{R}^{4n}$.

**Definition 5.5.** Let $(M, \psi)$ be a tetraplectic manifold on which the group $Q^k$ acts preserving $\psi$. We say that this action is *generalized Hamiltonian* if for any $X \in \mathfrak{q}$ the fundamental vector field $\hat{X}$ is Hamiltonian.

The following basis of $\mathfrak{q}$, representing the quaternion imaginary units $i$, $j$, $k$, is convenient for describing the tri-moment map;

$$H = \begin{pmatrix} i & 0 \\ 0 & -i \end{pmatrix}, \quad X = \begin{pmatrix} 0 & 1 \\ -1 & 0 \end{pmatrix}, \quad Y = \begin{pmatrix} 0 & i \\ i & 0 \end{pmatrix}.$$

The space of 3-vectors $\Lambda^3(\mathfrak{q})$ can be identified with $\mathbb{R}$ by the isomorphism that sends $X \wedge Y \wedge H \mapsto 1$.

Any

$$\delta = (\delta_1, \delta_2, \ldots, \delta_k) = (U_1 \wedge V_1 \wedge W_1, \ldots, U_k \wedge V_k \wedge W_k) \in (\Lambda^3 \mathfrak{q})^k$$

induces a $k$-tuple of 3-vector fields on $M$

$$e_\delta = (e_{\delta_1}, e_{\delta_2}, \ldots, e_{\delta_k}) = (\hat{U}_1 \wedge \hat{V}_1 \wedge \hat{W}_1, \ldots, \hat{U}_k \wedge \hat{V}_k \wedge \hat{W}_k) \in (\Lambda^3(TM))^k.$$

The following definition generalizes the classical moment map.

**Definition 5.6.** Let $Q^k$ act on a tetraplectic manifold $(M, \psi)$ in a generalized Hamiltonian fashion. A *tri-moment map* $\mu$ is a map

$$\mu : M \to (\Lambda^3 \mathfrak{q}^k)^* \cong (\Lambda^3 \mathfrak{q}^*)^k \cong \mathbb{R}^k$$

satisfying the following conditions:

1. $\mu$ is $Q^k$-invariant, i.e. $\mu(g \cdot p) = \mu(p)$. This is equivalent to saying that $\mu$ is $Q^k$-equivariant because the coadjoint action of $Q$ on $\mathfrak{q}^*$ induces a trivial action of $Q$ on $\Lambda^3 \mathfrak{q}^*$.

2. For any $\delta = (\delta_1, \ldots, \delta_k) \in (\Lambda^3 \mathfrak{q})^k$, $p \in M$ and $v \in T_p M$, we have

$$d\mu_p(v)(\delta) = \sum_{i=1}^{k} \iota_{e_{\delta_i}(p)} \psi(v) =: \iota_{e_\delta(p)} \psi(v),$$

where $e_{\delta_i}$ is the 3-vector field induced by $\delta_i$.

The following example is explained in detail in [Gentili et al., 2019].

**Example 5.7** (Standard model). Let $Q^n = \mathrm{Sp}(1)^n$ act on $\mathbb{H}^n$ by component-wise left multiplication (a.k.a the regular representation):

$$(q_1, \ldots, q_n) \cdot (h_1, \ldots, h_n) := (q_1 h_1, \ldots, q_n h_n).$$

Then this action is generalized Hamiltonian with respect to the standard tetraplectic form

$$\psi_{\mathbb{H}^n} := \sum_{i=1}^{n} dx_{4i-3} \wedge dx_{4i-2} \wedge dx_{4i-1} \wedge dx_{4i},$$

where $x_1, \ldots, x_{4n}$ are real coordinates on $\mathbb{H}^n \cong \mathbb{R}^{4n}$.

The associated tri-moment map is given by

$$\mu_{\mathbb{H}^n}(h_1, \ldots, h_m) = -\frac{1}{4}\left(|h_1|^4, \ldots, |h_n|^4\right) + C,$$

for some constant $C \in \mathbb{R}^n$.

**Remark 5.8.** This example provides the local model for generalized Hamiltonian $Q^n$-actions. The fibers of the tri-moment map are $Q^n$-orbits, and the image lies in $\mathbb{R}^n_{\geq}$, making this a quaternionic analogue of the standard symplectic toric model.

## 5.2 Generalized Lagrangian-type toric fibrations

The following definition is a "natural" generalization of a tri-moment map.

**Definition 5.9.** Let $(M^{4n}, \psi)$ be a a smooth $4n$-dimensional manifold with a tetraplectic structure, and let $B$ be an $n$-dimensional manifold with corners. A map $\pi : (M^{4n}, \psi) \to B$ is called a *locally generalized Lagrangian-type toric fibration* if there exists a system $\{(U_\alpha^B, \phi_\alpha^B)\}_{\alpha \in \mathcal{A}}$ of coordinate neighborhoods of $B$ modeled on $\mathbb{R}^n_{\geq}$ and for each index $\alpha$ there exists a tetraplectomorphism $\phi_\alpha^M : (\pi^{-1}(U_\alpha^B), \psi_M) \to (\pi_{\mathbb{H}^n}^{-1}(U_\alpha^B), \psi_{\mathbb{H}^n})$ such that $\pi_{\mathbb{H}^n} \circ \phi_\alpha^M = \phi_\alpha^B \circ \pi$.

**Example 5.10.** A tri-moment map associated with a quaternionic toric manifold, as introduced in [Gentili et al., 2019] defines a locally generalized Lagrangian-type toric fibration. In this setting, the manifold $M^{4n}$ is equipped with a tetraplectic structure $\psi$, and a $Q^n = \mathrm{Sp}(1)^n$-action that is generalized Hamiltonian. The corresponding tri-moment map

$$\mu : (M^{4n}, \psi) \to \mathbb{R}^n$$

satisfies the properties of a locally generalized Lagrangian-type toric fibration, as the fibers are $Q^n$-orbits, and locally, the structure is modeled on the regular $Q^n$-action on $\mathbb{H}^n$ with its canonical tetraplectic form.

**Example 5.11.** Let $Q^n = (S^3)^n$, and consider the product $\mathbb{R}^n \times Q^n$. Define the tetraplectic form on this space by

$$\psi_{\mathbb{R}^n \times Q^n} = \sum_{k=1}^n d\xi_k \wedge d\theta_k^1 \wedge d\theta_k^2 \wedge d\theta_k^3,$$

where $(\xi_1, \ldots, \xi_n)$ are coordinates on $\mathbb{R}^n$, and for each $k$, the $(\theta_k^1, \theta_k^2, \theta_k^3)$ are angular coordinates on the $k$-th copy of $S^3$.

Then the projection
$$\mathrm{pr}_1 : (\mathbb{R}^n \times Q^n, \psi_{\mathbb{R}^n \times Q^n}) \to \mathbb{R}^n$$
onto the first factor defines a locally generalized Lagrangian-type toric fibration. Here, the coordinate neighborhoods $U_\alpha^B \subset \mathbb{R}^n$ with identity maps $\phi_\alpha^B = \mathrm{id}$, and the corresponding tetraplectomorphisms $\phi_\alpha^M = \mathrm{id}$ give the trivializations required in the definition (under the identification in polar coordinates).

**Lemma 5.12** (Generalized Liouville Coordinates). *Let $\pi : (M^{4n}, \psi) \to B$ be a locally generalized Lagrangian-type toric fibration with fibers diffeomorphic to quaternionic tori $Q^n = (S^3)^n$. Suppose there exists a local Lagrangian section $f$ of $\mathrm{pr}_1 : (\mathbb{R}^n \times Q^n, \psi_{\mathbb{R}^n \times Q^n}) \to \mathbb{R}^n$, i.e. $f^*\psi_{\mathbb{R}^n \times Q^n} = 0$, where*

$$\psi_{\mathbb{R}^n \times Q^n} = \sum_{k=1}^n d\xi_k \wedge d\theta_k^1 \wedge d\theta_k^2 \wedge d\theta_k^3.$$

*Then the map*
$$\Psi : V \times Q^n \to \pi^{-1}(V), \quad \Psi(\xi, u) := (\xi, uf(\xi))$$
*defines a local tetraplectomorphism, and the coordinates $(\xi_k, \theta_k^1, \theta_k^2, \theta_k^3)$ are quaternionic analogues of action-angle coordinates.*

*Proof.* The map $\Psi$ is a diffeomorphism from $V \times Q^n$ onto itself (since it's smooth, invertible with smooth inverse $\Psi^{-1}(\xi, v) = (\xi, vf(\xi)^{-1})$). The fiber-preserving nature of $\Psi$ is immediate from its definition, and the group multiplication on $Q^n$ ensures that $uf(\xi) \in Q^n$.

To check that $\Psi$ preserves the tetraplectic structure, we compute the pullback:

$$\Psi^*\psi_{\mathbb{R}^n \times Q^n} = \sum_{k=1}^n d\xi_k \wedge d(\theta_k^1 + \theta_k^1(f(\xi))) \wedge d(\theta_k^2 + \theta_k^2(f(\xi))) \wedge d(\theta_k^3 + \theta_k^3(f(\xi))).$$

By the assumption $f^*\psi_{\mathbb{R}^n \times Q^n} = 0$, the terms involving derivatives of $f$ cancel out in the wedge product terms, and so
$$\Psi^*\psi_{\mathbb{R}^n \times Q^n} = \psi_{\mathbb{R}^n \times Q^n}.$$

Hence, $\Psi$ is a tetraplectomorphism. □

**Remark 5.13.** This result is a quaternionic–tetraplectic analogue of the classical Liouville theorem in integrable Hamiltonian systems. In symplectic geometry, the image of a Lagrangian section defines local action-angle coordinates via the identification of the cotangent bundle with a product of coordinates. Here, the section $f$ is called Lagrangian because it satisfies $f^*\psi = 0$, and the resulting trivialization $\Psi$ of the bundle is a local tetraplectomorphism.

The following Generalized Arnold-Liouville theorem says that $\mathrm{pr}_1 : (\mathbb{R}^n \times Q^n, \psi_{\mathbb{R}^n \times Q^n}) \to \mathbb{R}^n$ is the local model of a locally generalized Lagrangian-type toric fibration.

**Theorem.** *(Generalized Arnold-Liouville theorem.) Let $F$ be a fiber of a locally generalized Lagrangian-type toric fibration $\pi : (M^{4n}, \psi) \to B$. Then there exists an open neighborhood of $F$ that is fiber-preserving tetraplectomorphic to $(V \times Q^n, \psi_{\mathbb{R}^n \times Q^n}) \to \mathbb{R}^n$, where $V \subset \mathbb{R}^n$.*

*Proof.* Let $\pi : (M^{4n}, \psi) \to B$ be a locally generalized Lagrangian-type toric fibration, and let $F = \pi^{-1}(b)$ be a fiber over a point $b \in B$. By Definition 2.3, there exists an open neighborhood $U_b^B \subset B$ of $b$, together with a coordinate chart $\varphi_b^B : U_b^B \to \mathbb{R}^n$ and a tetraplectomorphism

$$\varphi_b^M : (\pi^{-1}(U_b^B), \psi) \to (\pi_{\mathbb{H}^n}^{-1}(\varphi_b^B(U_b^B)), \psi_{\mathbb{H}^n}),$$

such that the following diagram commutes:

$$\begin{array}{ccc} \pi^{-1}(U_b^B) & \xrightarrow{\varphi_b^M} & \pi_{\mathbb{H}^n}^{-1}(\varphi_b^B(U_b^B)) \\ \pi \downarrow & & \downarrow \pi_{\mathbb{H}^n} \\ U_b^B & \xrightarrow{\varphi_b^B} & \phi^B(U_b^B) \end{array}$$

By construction of the standard model 3.8, we identify $\pi_{\mathbb{H}^n}^{-1}(\phi^B(U_b^B)) \cong \phi^B(U_b^B) \times Q^n$ (see also [Hopkinson, 2012, Section 3.2] for details on the identification), with tetraplectic structure $\psi_{\mathbb{R}^n \times Q^n}$ from Example 5.11. The map $\varphi_b^M$ is a diffeomorphism that preserves the 4-form $\psi$ and intertwines the fibrations, i.e., $\pi_{\mathbb{H}^n} \circ \varphi_b^M = \varphi_b^B \circ \pi$.

Now take $V := \varphi_b^B(U_b^B) \subset \mathbb{R}^n$. Then we have that

$$\varphi_b^M : \pi^{-1}(U_b^B) \to V \times Q^n$$

is a tetraplectomorphism over $V$, i.e.,

$$(\varphi_b^M)^* \psi_{\mathbb{R}^n \times Q^n} = \psi.$$

Therefore, the neighborhood $\pi^{-1}(U_b^B)$ of the fiber $F$ is fiber-preserving tetraplectomorphic to $V \times Q^n$, which completes the proof. $\square$

**Remark 5.14.** As a consequence of the Generalized Arnold–Liouville theorem, notice that the fibers of a locally generalized Lagrangian-type toric fibration are diffeomorphic to quaternionic tori $Q^n = (S^3)^n$. These fibers generalize the classical notion of Lagrangian tori in symplectic geometry.

**Remark 5.15** (Comparison with the classical Arnold–Liouville theorem)**.** The classical Arnold–Liouville theorem in symplectic geometry asserts that a completely integrable system on a $2n$-dimensional symplectic manifold admits a local model as a Lagrangian torus fibration over $\mathbb{R}^n$, with fibers diffeomorphic to $T^n$ and a canonical expression for the symplectic form in action-angle coordinates.

In our quaternionic-tetraplectic setting, we replace the symplectic 2-form by a closed, non-degenerate 4-form $\psi$, and the torus $T^n$ by the quaternionic torus $Q^n = \text{Sp}(1)^n \cong (S^3)^n$. The notion of a completely integrable system is generalized by the existence of a tri-moment map associated with a local $Q^n$-action preserving the tetraplectic form. The fibers of this fibration are $Q^n$, and the local model becomes

$$\text{pr}_1 : (\mathbb{R}^n \times Q^n, \psi_{\mathbb{R}^n \times Q^n}) \to \mathbb{R}^n.$$

Thus, the Generalized Arnold–Liouville theorem established here plays an analogous role, providing a local normal form for tetraplectic manifolds with locally generalized Lagrangian-type toric fibrations. However, classical action-angle coordinates do not directly generalize due to the noncommutativity of $Q^n$.

**Definition 5.16** (Fiber-preserving tetraplectomorphism)**.** Let $\pi_1 : (M_1^{4n}, \psi_1) \to B$ and $\pi_2 : (M_2^{4n}, \psi_2) \to B$ be fibrations of smooth manifolds equipped with tetraplectic forms $\psi_1, \psi_2$. A diffeomorphism

$$f : M_1^{4n} \to M_2^{4n}$$

is called a *fiber-preserving tetraplectomorphism* if:

a) $f^*\psi_2 = \psi_1$, i.e., $f$ is a tetraplectomorphism,

b) $\pi_2 \circ f = \pi_1$, i.e., $f$ preserves (diffeomorphically) the fibers over $B$.

As an application, a diffeomorphism

$$\phi : \mathbb{R}^n \times Q^n \to \mathbb{R}^n \times Q^n$$

we will be called *fiber-preserving tetraplectomorphism of the standard model* applying the above definition to $\psi_1 = \psi_2 = \psi_{\mathbb{R}^n \times Q^n}$, and $\pi_1 = \pi_2 = \mathrm{pr}_1 : (\mathbb{R}^n \times Q^n, \psi_{\mathbb{R}^n \times Q^n}) \to \mathbb{R}^n$.

**Lemma 5.17.** *Let $\phi$ be a fiber-preserving tetraplectomorphism of the standard model. Then $\phi$ is of the form*

$$\phi(\xi, u) = (\rho^{-T}(\xi) + c, \rho(u)f(\xi))$$

*for some $\rho \in \mathrm{Aut}(Q^n)$, $\rho^{-T}$ the transpose inverse of $\rho$, $c \in \mathbb{R}^n$, and a Lagrangian section $f$ of $\mathrm{pr}_1 : (\mathbb{R}^n \times Q^n, \psi_{\mathbb{R}^n \times Q^n}) \to \mathbb{R}^n$.*

*Proof.* Suppose $\phi : \mathbb{R}^n \times Q^n \to \mathbb{R}^n \times Q^n$ is a fiber-preserving diffeomorphism such that $\phi^*\psi_{\mathbb{R}^n \times Q^n} = \psi_{\mathbb{R}^n \times Q^n}$. Since $\phi$ is fiber-preserving, we can write

$$\phi(\xi, u) = (\phi_B(\xi), \phi_F(\xi, u)),$$

where $\phi_B : \mathbb{R}^n \to \mathbb{R}^n$ and for each fixed $\xi$, the map $u \mapsto \phi_F(\xi, u)$ is a diffeomorphism of $Q^n$ onto the fiber over $\phi_B(\xi)$.

Because $Q^n$ is a Lie group, any such fiber map must be of the form

$$\phi_F(\xi, u) = \rho(u)f(\xi)$$

for some $\rho \in \mathrm{Aut}(Q^n)$ and some smooth function $f : \mathbb{R}^n \to \mathbb{R}^n \times Q^n$. Therefore, $\phi$ has the form

$$\phi(\xi, u) = (\phi_B(\xi), \rho(u)f(\xi)).$$

Now, we impose the tetraplectomorphism condition: $\phi^*\psi_{\mathbb{R}^n \times Q^n} = \psi_{\mathbb{R}^n \times Q^n}$. Under $\phi$, the pullback of $d\xi_k$ becomes a linear combination of $d\xi_j$ via the Jacobian of $\phi_B$, and the $d\theta_k^j$ are left-multiplied by $\rho$ and twisted by $f(\xi)$.

To preserve $\psi$, we must have:

- $\phi_B(\xi) = \rho^{-T}(\xi) + c$ for some $c \in \mathbb{R}^n$, so that the pullback of $d\xi_k$ matches the transformation law under $\rho$,

- $f^*\psi_{\mathbb{R}^n \times Q^n} = 0$.

This implies that $f$ is a *Lagrangian section* with respect to the tetraplectic form. Therefore, any fiber-preserving tetraplectomorphism has the desired form. $\square$

**Proposition 5.18.** *Let $\pi : (M^{4n}, \psi) \to B^n$ be a locally generalized Lagrangian-type fibration on an $n$-dimensional base $B^n$, and let $\{(U_\alpha^B, \varphi_\alpha^B, \varphi_\alpha^M)\}$ be a system of local identifications of $\pi$ with the standard model $\pi_{\mathbb{H}^n} : (\mathbb{H}^n, \psi_{\mathbb{H}^n}) \to \mathbb{R}_{\geq 0}^n$. Then, on each nonempty overlap $U_{\alpha\beta}^B := U_\alpha^B \cap U_\beta^B$, there exists an automorphism $\rho_{\alpha\beta} \in \mathrm{Aut}(Q^n)$ such that the overlap map*

$$\varphi_\alpha^M \circ (\varphi_\beta^M)^{-1} : \pi_{\mathbb{H}^n}^{-1}(\varphi_\beta^B(U_{\alpha\beta}^B)) \to \pi_{\mathbb{H}^n}^{-1}(\varphi_\alpha^B(U_{\alpha\beta}^B))$$

*is $\rho_{\alpha\beta}$-equivariant.*

*Proof.* We first focus on the interior $B^\circ := B \setminus \partial B$. Since the restriction of $\pi$ to $B^\circ$ is a (nonsingular) locally generalized Lagrangian-type toric fibration, it is locally identified with the standard model.

For each $\alpha$, the chart $\varphi_\alpha^M$ gives a fiber-preserving tetraplectomorphism

$$\varphi_\alpha^M : (\pi^{-1}(U_\alpha^B \setminus \partial B), \psi) \xrightarrow{\cong} (\pi_{\mathbb{H}^n}^{-1}(\varphi_\alpha^B(U_\alpha^B \setminus \partial B)), \psi_{\mathbb{H}^n}),$$

which covers the diffeomorphism $\varphi_\alpha^B : U_\alpha^B \to \varphi_\alpha^B(U_\alpha^B) \subset \mathbb{R}^n_{\geq}$.

We define a standard model identification

$$\nu_{\mathbb{H}^n} : \varphi_\alpha^B(U_\alpha^B \setminus \partial B) \times Q^n \to \pi_{\mathbb{H}^n}^{-1}(\varphi_\alpha^B(U_\alpha^B \setminus \partial B))$$

by

$$\nu_{\mathbb{H}^n}(\xi, u) := u \cdot \iota(\xi),$$

where $\iota : \mathbb{R}^n_{\geq} \to \mathbb{H}^n$ is the section defined by $\iota(\xi_1, \ldots, \xi_n) := (\xi_1^{1/4}, \ldots, \xi_n^{1/4}) \in \mathbb{H}^n$. This map is $Q^n$-equivariant with respect to the multiplicative action on the second factor. Hence, $\nu_{\mathbb{H}^n}^{-1} \circ \varphi_\alpha^M$ is a local identification.

Now, consider a nonempty overlap $U_{\alpha\beta}^B$. By applying the same technique as above, we compare the charts via:

$$\nu_{\mathbb{H}^n}^{-1} \circ \varphi_\alpha^M \circ (\varphi_\beta^M)^{-1} \circ \nu_{\mathbb{H}^n} : \varphi_\beta^B(U_{\alpha\beta}^B) \times Q^n \to \varphi_\alpha^B(U_{\alpha\beta}^B) \times Q^n.$$

Since both $\varphi_\alpha^M$ and $\varphi_\beta^M$ are fiber-preserving tetraplectomorphisms and $\nu_{\mathbb{H}^n}$ is $Q^n$-equivariant, the transition map is a $Q^n$-equivariant diffeomorphism of the form

$$(\xi, u) \mapsto (\rho_{\alpha\beta}(\xi), \rho(u))$$

for some $\rho \in \mathrm{Aut}(Q^n)$. Thus, we conclude that the overlap map

$$\varphi_\alpha^M \circ (\varphi_\beta^M)^{-1}$$

is $\rho_{\alpha\beta}$-equivariant on $\pi^{-1}(U_{\alpha\beta})$, and $\rho_{\alpha\beta} \in \mathrm{Aut}(Q^n)$. □

## 5.3 Integral affine structures

Definition 2.3 introduces a rigid structure for manifolds with local $Q^n$-actions, where the local charts are required to be compatible up to $S^3$-symmetries. This mirrors the complex toric case, where local charts are compatible up to $T^n$-symmetries and transition maps lie in the group of integral affine transformations.

In that setting, the transition maps form the affine group $GL(n, \mathbb{Z}) \ltimes \mathbb{R}^n$, with the linear part capturing symmetries of the torus and the translational part accounting for changes in base coordinates. The analogue in the quaternionic case is given by the group of automorphisms of the quaternionic torus, $\mathrm{Aut}(Q^n)$, which serves as the linear part of the structure group.

Here, to describe transitions between charts, we consider the *quaternionic affine group*:

$$Aff_Q(n) := Aut(Q^n) \ltimes \mathbb{H}^n,$$

where $Aut(Q^n)$ corresponds to the linear part, while $\mathbb{H}^n$ captures translations. If we restrict attention to integral structures (analogous to the lattice $\mathbb{Z}^n$ in the toric case), the discrete model becomes

$$GL(n, 4\pi\mathbb{Z}) \ltimes (4\pi\mathbb{Z})^n,$$

since the exponential map for $SU(2)$ has kernel generated by $4\pi$-rotations.

We now formalize this structure:

**Definition 5.19.** A *quaternionic integral affine structure* on an $n$-dimensional manifold with corners $B$ is an atlas $\{(U_\alpha^B, \varphi_\alpha^B)\}$ whose transition maps $\varphi_{\alpha\beta}^B := \varphi_\alpha^B \circ (\varphi_\beta^B)^{-1}$ lie in $Aff_Q(n)$, i.e., each overlap map belongs to $Aut(Q^n) \ltimes \mathbb{H}^n$.

**Lemma 5.20.** *Let $B$ be an $n$-dimensional $C^\infty$ manifold with corners equipped with a quaternionic integral affine structure $\{(U_\alpha^B, \phi_\alpha^B)\}$. Then there exists a characteristic pair $(P, \mathcal{L})$ on $B$ associated with this atlas.*

*Proof.* By definition, a quaternionic integral affine structure means that the transition functions $\phi_{\alpha\beta}^B := \phi_\alpha^B \circ (\phi_\beta^B)^{-1}$ take values in the affine group $\mathrm{Aut}(Q^n) \ltimes (4\pi\mathbb{Z})^n$, where $\mathrm{Aut}(Q^n) \cong S_n \ltimes (\mathbb{Z}/2\mathbb{Z})^n$ acts discretely.

This reduction of the structure group determines a principal $\mathrm{Aut}(Q^n)$-bundle $P \to B$, which plays the role of the quaternionic frame bundle. Associated to this is a rank-$n$ lattice bundle $\pi_\Lambda : \Lambda_P \to B$ whose fibers are isomorphic to the quaternionic integral lattice $\Lambda_Q := (4\pi\mathbb{Z})^n \subset \mathbb{R}^n$.

Let $\mathcal{S}^{(n-1)}B$ denote the codimension-one strata of $B$. On a coordinate neighborhood $(U_\alpha^B, \phi_\alpha^B)$ such that $U_\alpha^B \cap \mathcal{S}^{(n-1)}B \neq \emptyset$, the image $\phi_\alpha^B(U_\alpha^B \cap \mathcal{S}^{(n-1)}B)$ lies in a hyperplane $\langle \xi, u_\alpha \rangle = 0$ for a primitive vector $u_\alpha \in \Lambda_Q$.

For another chart $(U_\beta^B, \phi_\beta^B)$ intersecting the same face, the corresponding primitive vector $u_\beta \in \Lambda_Q$ satisfies $u_\alpha = \rho_{\alpha\beta}^{-T}(u_\beta)$, where $\rho_{\alpha\beta} \in \mathrm{Aut}(Q^n)$ is the linear part of the transition map. Since the transition automorphisms preserve the lattice structure, $u_\alpha$ and $u_\beta$ generate the same sublattice of $\Lambda_Q$.

This allows us to define, over $\mathcal{S}^{(n-1)}B$, a rank-one sublattice bundle $\mathcal{L} \to \mathcal{S}^{(n-1)}B$, locally generated by $u_\alpha$, which is patched compatibly under $\mathrm{Aut}(Q^n)$. By construction, the data of $(P, \mathcal{L})$ satisfies the unimodularity condition: over each $k$-dimensional stratum $\mathcal{S}^{(k)}B$, the corresponding $n-k$ primitive vectors span a direct summand of rank $n-k$ in $\Lambda_Q$.

Therefore, the pair $(P, \mathcal{L})$ is a characteristic pair on $B$ associated to the quaternionic integral affine atlas $\{(U_\alpha^B, \phi_\alpha^B)\}$. $\square$

**Lemma 5.21.** *The canonical model $M_{(P,\mathcal{L})}$ is a $C^\infty$ manifold and admits a tetraplectic structure $\psi$ such that the projection $\pi : M_{(P,\mathcal{L})} \to B$ is a locally generalized Lagrangian-type toric fibration.*

*Proof.* We reconstruct $M_{(P,\mathcal{L})}$ from the total space $Q_P$ of the principal $Q^n$-bundle $P \to B$, using the tetraplectic cutting technique described in [Gentili et al., 2019, Proposition 5.2], analogous to the symplectic cut in [Yoshida, 2011, Lemma 7.3].

Let $\{(U_\alpha^B, \varphi_\alpha^B)\}_{\alpha \in \mathcal{A}}$ be an integral affine atlas of $B$, and let $\varphi_\alpha^P : \pi_P^{-1}(U_\alpha^B) \cong U_\alpha^B \times Q^n$ be a local trivialization of $Q_P$. Each chart corresponds locally to $\mathbb{H}^n \cong \mathbb{R}^{4n}$, with standard coordinates $(\xi_1, \ldots, \xi_n) \in \mathbb{R}_{\geq 0}^n$ and quaternions $q_i \in \mathbb{H}$. The tetraplectic form is locally given by

$$\psi = \sum_{i=1}^n dx_{4i-3} \wedge dx_{4i-2} \wedge dx_{4i-1} \wedge dx_{4i},$$

where each $q_i \in \mathbb{H}$ is expressed in real coordinates as $(x_{4i-3}, \ldots, x_{4i})$.

For each $\alpha$ define the index set

$$I_\alpha := \{i \in \{1, \ldots, n\} \mid \varphi_\alpha^B(U_\alpha^B) \cap \{\xi_i = 0\} \neq \emptyset\},$$

i.e., those coordinate directions that vanish along a face of the polytope. For each $i \in I_\alpha$, define the subgroup $Z_i \subset Q^n$ acting only in the $i$-th quaternionic coordinate.

Define the subgroup $T_\alpha := \{u \in Q^n \mid u_i = 1 \text{ for } i \notin I_\alpha\} \subset Q^n$. It acts on $V_\alpha \times Q^n \times \mathbb{H}^{I_\alpha}$ by:

$$u \cdot (\xi, q, z) = (\xi, uq, (u_i^{-1} z_i)_{i \in I_\alpha}).$$

This action is generalized Hamiltonian with respect to the tetraplectic form

$$\psi_\alpha := \psi_{\mathbb{R}^n \times Q^n} + \sum_{i \in I_\alpha} d\bar{z}_i \wedge dz_i \wedge d\bar{z}_i \wedge dz_i,$$

and has a moment map

$$\Phi_\alpha(\xi, q, z) = (\xi_i - |z_i|^4)_{i \in I_\alpha}.$$

Let $\Phi_\alpha^{-1}(0)/T_\alpha$ denote the reduced space. Then the quotient space is a smooth manifold equipped with a tetraplectic form $\psi_\alpha^{\text{red}}$, and the map

$$\varphi_\alpha^M : \mu_{(P,\mathcal{L})}^{-1}(U_\alpha^B) \to \Phi_\alpha^{-1}(0)/T_\alpha$$

defines a local tetraplectomorphism.

Repeat this construction over all $\alpha \in \mathcal{A}$, and glue the charts $\Phi_\alpha^{-1}(0)/T_\alpha$ using the transition functions of $P$ (which lie in $\text{Aut}(Q^n)$), ensuring compatibility of tetraplectic forms.

This construction produces a global $C^\infty$ manifold $M_{(P,\mathcal{L})}$ equipped with a closed, non-degenerate 4-form $\psi$, invariant under the local $Q^n$-action. The projection

$$\pi : M_{(P,\mathcal{L})} \to B$$

is locally modeled on

$$\pi_{\mathbb{H}^n} : \mathbb{H}^n \to \mathbb{R}_{\geq 0}^n, \quad q \mapsto (|q_1|^4, \ldots, |q_n|^4),$$

and hence is a locally generalized Lagrangian-type toric fibration. □

**Remark 5.22.** The local tetraplectic cut used in the above proof is a quaternionic analogue of Lerman's symplectic cut, and is justified by the tetraplectic reduction procedure established by Gentili–Gori–Sarfatti [Gentili et al., 2019, Theorem 5.3]. In particular, the local charts of $M_{(P,\mathcal{L})}$ are constructed as reduced spaces of $\mathbb{H}^n \times \mathbb{H}^{|I_\alpha|}$ with respect to suitable $Q$-actions, as in their tetraplectic cut.

## 5.4 The necessary and sufficient condition

In this subsection we study when a 4n-dimensional manifold equipped with a local quaternionic torus action becomes a locally generalized Lagrangian-type toric fibration.

**Definition 5.23.** Let $(M^{4n}, \mathcal{Q})$ be a $4n$-dimensional manifold equipped with a local $Q^n$-action $\mathcal{Q}$. Let $\{(U_\alpha^M, \varphi_\alpha^M)\}_{\alpha \in \mathcal{A}}$ in $\mathcal{Q}$ be a weakly regular atlas of $M$, and let $\{(U_\alpha^B, \varphi_\alpha^B)\}_{\alpha \in \mathcal{A}}$ be the induced atlas of the orbit space $B_M$. For each nonempty overlap $U_{\alpha\beta}^M := U_\alpha^M \cap U_\beta^M$, let $\rho_{\alpha\beta} \in \text{Aut}(Q^n)$ be the transition automorphism of Definition 2.3. We say that the atlas $\{(U_\alpha^B, \varphi_\alpha^B)\}_{\alpha \in \mathcal{A}}$ defines a *quaternionic integral affine structure compatible with* $\{(U_\alpha^M, \varphi_\alpha^M)\}$ if, for each nonempty overlap $U_{\alpha\beta}^B := U_\alpha^B \cap U_\beta^B$, the linear part of the coordinate transition

$$A_{\alpha\beta} := D(\varphi_\alpha^B \circ (\varphi_\beta^B)^{-1})$$

belongs to $\text{GL}_n(4\pi\mathbb{Z})$ and satisfies

$$A_{\alpha\beta} = \rho_{\alpha\beta}^{-T}.$$

The next result is a consequence of Lemma 5.17 and Proposition 5.18.

**Lemma 5.24.** *If there exists a tetraplectic structure $\psi$ on $M^{4n}$ and a weakly regular atlas $\{(U_\alpha^M, \varphi_\alpha^M)\}_{\alpha \in \mathcal{A}} \in \mathcal{Q}$ such that on each $U_\alpha^M$, $\varphi_\alpha^M$ preserves tetraplectic structures, namely $\psi$, then the atlas $\{(U_\alpha^B, \varphi_\alpha^B)\}_{\alpha \in \mathcal{A}}$ of $B_M$ induced by $\{(U_\alpha^M, \varphi_\alpha^M)\}_{\alpha \in \mathcal{A}}$ is a quaternionic integral affine structure that is compatible with $\{(U_\alpha^M, \varphi_\alpha^M)\}_{\alpha \in \mathcal{A}}$. In particular, $B_M$ becomes a manifold with corners.*

Finally, summarizing the previous results, we get the following tetraplectic realization theorem.

**Theorem.** *Let $(M^{4n}, \mathcal{Q})$ be a $4n$-dimensional manifold equipped with a local $Q^n$-action $\mathcal{Q}$. There exists a tetrapletcic structure $\psi$ on $M^{4n}$ and there also exists a weakly regular atlas $\{(U_\alpha^M, \varphi_\alpha^M)\}_{\alpha \in \mathcal{A}} \in \mathcal{Q}$ of $M^{4n}$ such that on each $U_\alpha^M$, $\psi = (\varphi_\alpha^M)^* \psi_{\mathbb{H}^n}$ if and only if the atlas $\{(U_\alpha^B, \varphi_\alpha^B)\}_{\alpha \in \mathcal{A}}$ of $B_M$ induced by $\{(U_\alpha^M, \varphi_\alpha^M)\}_{\alpha \in \mathcal{A}}$ is a quaternionic integral affine structure compatible with $\{(U_\alpha^M, \varphi_\alpha^M)\}_{\alpha \in \mathcal{A}}$ and on each nonempty overlap $U_{\alpha\beta}^B$, the section $\theta_{\alpha\beta}^M$ is a Lagrangian, i.e. $(\theta_{\alpha\beta}^M)^* \psi_{Q_M}$ vanishes, where $\psi_{Q_M}$ is the tetraplectic structure defined in Lemma 5.21.*

## 5.5  Quaternionic Delzant-type Classification

Suppose that $(M^{4n}, \mathcal{Q})$ is a smooth compact $4n$-dimensional manifold equipped with a local $Q^n$-action $\mathcal{Q}$, such that the orbit space $B_M$ is a quaternionic integral affine manifold with corners, and the local "gluing" maps define Lagrangian overlaps, as in Theorem 5.4.

Then $M^{4n}$ admits a tetraplectic structure $\psi$ such that the orbit map

$$\pi_M : (M^{4n}, \psi) \to B_M$$

is a locally generalized Lagrangian-type toric fibration. The overlaps define a Čech cohomology class

$$\lambda(M) \in H^1(B_M, \mathscr{L}_{Q_M}^{\text{Lag}}),$$

called the *Lagrangian class*, where $\mathscr{L}_{Q_M}^{\text{Lag}}$ is the sheaf of Lagrangian sections of the torus bundle $\pi_{Q_M} : (Q_M, \psi_{Q_M}) \to B_M$.

**Theorem** (Quaternionic Delzant-type Classification). *Let $(M_i^{4n}, \mathcal{Q}_i, \psi_i)$, $i = 1, 2$, be compact smooth $4n$-manifolds equipped with local $Q^n$-actions and tetraplectic structures, such that $\pi_i : (M_i^{4n}, \psi_i) \to B_{M_i}$ are locally generalized Lagrangian-type toric fibrations. Then $M_1$ and $M_2$ are fiber-preserving tetraplectomorphic if and only if there exists a diffeomorphism $f_B : B_{M_1} \to B_{M_2}$ which preserves the quaternionic integral affine structures and satisfies $f_B^* \lambda(M_2) = \lambda(M_1)$, where $\lambda(M_i) \in H^1(B_{M_i}, \mathscr{L}_{Q_{M_i}}^{\text{Lag}})$ are the respective Lagrangian classes.*

*Moreover, for any compact quaternionic integral affine manifold $B$ with corners and any class $\lambda \in H^1(B, \mathscr{L}_{Q_B}^{\text{Lag}})$, there exists a locally generalized Lagrangian-type fibration $\pi : (M^{4n}, \psi) \to B$ realizing $\lambda$ as its Lagrangian class. The pair $(B, \lambda)$ classifies $M^{4n}$ up to fiber-preserving tetraplectomorphism.*

*Proof.* We split the proof into two parts.

**(1) Uniqueness.** Suppose there exists a diffeomorphism $f_B : B_{M_1} \to B_{M_2}$ preserving the quaternionic integral affine structures, and such that the pullback of the Lagrangian class satisfies

$$f_B^* \lambda(M_2) = \lambda(M_1).$$

Since the quaternionic affine structures agree under $f_B$, the associated $Q^n$-bundles $Q_{M_1}$ and $Q_{M_2}$ are isomorphic as topological group bundles.

Moreover, the Čech cocycles defining the gluing of local models in both $M_1$ and $M_2$ are identified under the pullback by $f_B$, because the Lagrangian class is preserved. As a result, the canonical models constructed from these data are diffeomorphic via a fiber-preserving tetraplectomorphism. Thus, $M_1 \cong M_2$ as tetraplectic manifolds with locally regular $Q^n$-actions.

**(2) Existence.** Let $B$ be a quaternionic integral affine manifold with corners, and $\lambda \in H^1(B, \mathscr{L}_Q^{\text{Lag}})$ a cohomology class. The affine structure allows us to construct a principal $Q^n$-bundle $Q_B \to B$ with local trivializations modeled on $\mathbb{R}^n \times Q^n \to \mathbb{R}^n$, and each equipped with the standard tetraplectic form

$$\psi_{\mathbb{R}^n \times Q^n} = \sum_{k=1}^{n} d\xi_k \wedge d\theta_k^1 \wedge d\theta_k^2 \wedge d\theta_k^3.$$

The Lagrangian class $\lambda$ determines Čech cocycles that define the gluing of these charts via fiber-preserving tetraplectomorphisms. The vanishing of the pullback of $\psi$ along the transition sections guarantees that the global form is well-defined and smooth.

Thus, the resulting manifold $M^{4n}$ admits a locally regular $Q^n$-action, a tetraplectic form $\psi$, and a projection $\pi : M \to B$ which is a locally generalized Lagrangian-type toric fibration with Lagrangian class $\lambda$. $\square$

**Corollary 5.25.** *Two locally generalized Lagrangian-type toric fibrations $(M_1, \psi_1), (M_2, \psi_2)$ are fiber-preserving tetraplectomorphic if and only if their bases admit a diffeomorphism preserving the quaternionic integral affine structures and matching their Lagrangian classes.*

**Remark 5.26** (Comparison with the classical Delzant theorem). In the symplectic setting, the Delzant theorem establishes a one-to-one correspondence between compact symplectic toric manifolds and Delzant polytopes: simple, rational, smooth polytopes in $\mathbb{R}^n$ with a compatible integral affine structure.

Our quaternionic-tetraplectic classification plays a parallel role: it classifies compact $4n$-manifolds with local $Q^n$-actions and tetraplectic structures via two types of combinatorial-topological data:

1. A *quaternionic integral affine manifold $B$* with corners, generalizing the polytope;

2. A *Lagrangian class $\lambda \in H^1(B, \mathscr{L}_Q^{\mathrm{Lag}})$*, encoding the "gluing" of local models.

# References

[Atiyah, 1982] Atiyah, M. (1982). Convexity and commuting hamiltonians. *Bulletin of the London Mathematical Society*, 14:1–15.

[Buchstaber and Panov, 2000] Buchstaber, V. M. and Panov, T. E. (2000). Torus actions, combinatorial topology, and homological algebra. *Uspekhi Matematicheskikh Nauk*, 55(5(335)):3–106. Russian Math. Surveys 55:5 (2000), 825–921.

[Buchstaber and Panov, 2002] Buchstaber, V. M. and Panov, T. E. (2002). *Torus actions and their applications in topology and combinatorics*, volume 24 of *University Lecture Series*. American Mathematical Society, Providence, RI.

[Buchstaber and Panov, 2014] Buchstaber, V. M. and Panov, T. E. (2014). Toric topology.

[Davis, 2008] Davis, M. W. (2008). *The geometry and topology of Coxeter groups*, volume 32 of *London Mathematical Society Monographs Series*. Princeton University Press, Princeton, NJ.

[Davis and Januszkiewicz, 1991] Davis, M. W. and Januszkiewicz, T. (1991). Convex polytopes, coxeter orbifolds and torus actions. *Duke Mathematical Journal*, 62(2):417–451.

[Delzant, 1988] Delzant, T. (1988). Hamiltoniens périodiques et images convexes de l'application moment. *Bulletin de la Société Mathématique de France*, 116:315–339.

[Duistermaat, 1980] Duistermaat, J. J. (1980). On global action-angle coordinates. *Communications on Pure and Applied Mathematics*, 33(6):687–706.

[Foth, 2002] Foth, P. (2002). Tetraplectic structures, tri-momentum maps, and quaternionic flag manifolds. *Journal of Geometry and Physics*, 41:330–343.

[Gay and Symington, 2009] Gay, D. and Symington, M. (2009). Toric structures on near-symplectic 4-manifolds. *Journal of the European Mathematical Society (JEMS)*, 11:487–520.

[Gentili et al., 2019] Gentili, G., Gori, A., and Sarfatti, G. (2019). Quaternionic toric manifolds. *Journal of Symplectic Geometry*, 17(1):267–300.

[Gkeneralis and Prassidis, 2025] Gkeneralis, I. and Prassidis, S. (2025). Topological rigidity of quoric manifolds. *Colloquium Mathematicum*. To appear.

[Guillemin and Sternberg, 1982] Guillemin, V. and Sternberg, S. (1982). Convexity properties of the moment mapping. *Inventiones Mathematicae*, 67:491–513.

[Hatcher, 2017] Hatcher, A. (2017). *Vector Bundles and K-Theory (version 2.2)*. Available at https://pi.math.cornell.edu/~hatcher/VBKT/VBpage.html.

[Hopkinson, 2012] Hopkinson, J. (2012). *Quoric Manifolds*. PhD thesis, University of Manchester.

[Lawson and Michelsohn, 1989] Lawson, H. and Michelsohn, M. (1989). *Spin geometry*. Princeton University Press.

[Mishachev, 1996] Mishachev, K. (1996). The classification of lagrangian bundles over surface. *Differential Geometry and its Applications*, 6:301–320.

[Scott, 1995] Scott, R. (1995). Quaternionic toric varieties. *Duke Mathematical Journal*, 78(2):373–397.

[Symington, 2003] Symington, M. (2003). Four dimensions from two in symplectic topology. In *Topology and Geometry of Manifolds*, volume 71 of *Proc. Sympos. Pure Math.*, pages 153–208. Amer. Math. Soc., Providence, RI. Athens, GA, 2001.

[Yoshida, 2011] Yoshida, T. (2011). Local torus actions modeled on the standard representation. *Advances in Mathematics*, 227:1914–1955.